\documentclass[leqno]{article}
\usepackage{latexsym,amsmath,amsfonts,mathrsfs,amssymb,amsthm}
\usepackage{yhmath}
\usepackage[usenames]{color}
\usepackage{graphicx}
\def\R{\mathbb R}
\def\N{\mathbb N}
\def\mS{\mathbb S}
\def\C{{\cal C}}
\def\Z{\mathbb Z}
\def\s{\sharp}
\def\a{\alpha}
\def\b{\beta}
\def\e{\epsilon}
\def\d{\delta}
\def\t{\tilde}

\def\O{\Omega}
\def\cO{{\cal O}}
\def\mo{{\mathfrak o}}
\def\Y{\mathbb Y}
\def\T{\mathbb T}
\def\P{\mathbb P}
\def\H{{\cal H}}
\def\cV{{\cal V}}
\def\cK{{\cal K}}
\def\cC{{\cal C}}
\def\cD{{\cal D}}
\def\cT{{\cal T}}
\def\cU{{\cal U}}
\def\cG{{\cal G}}
\def\cI{{\cal I}}
\def\F{\mathfrak F}
\def\be{\begin{equation}}
\def\ee{\end{equation}}
\def\bs{\backslash}
\def\pa{\partial}
\def\rg{\rangle}
\def\lg{\langle}
\def\lf{\lfloor}
\def\qed{\hfill$\Box$\bigskip}
\def\nd{\noindent\textbf{Proof.} }
\def\tb{\textcolor{blue}}

\def\wg{\wedge}
\def\spt{\mbox{spt}}
\numberwithin{equation}{section}
\newtheorem{thm}{Theorem}[section]

\newtheorem{lem}[thm]{Lemma}
\newtheorem{pro}[thm]{Proposition}
\newtheorem{defn}[thm]{Definition}
\newtheorem{exam}[thm]{Example}
\newtheorem{cor}[thm]{Corollary}
\newtheorem{rem}[thm]{Remark}
\begin{document}
\bigskip

\begin{center}{\Large \textbf{On the Almgren minimality of the product of a paired calibrated set with a calibrated set of codimension 1 with singularities, and new Almgren minimal cones}}\end{center}

\centerline{\large Xiangyu Liang}
\centerline{maizeliang@gmail.com}
\centerline{School of mathematical sciences, Beihang University, P.R. China}

\bigskip

\centerline {\large\textbf{Abstract.}} 

In this paper, we prove that the product of a paired calibrated set and a set of codimension 1 calibrated by a coflat calibration with small singularity set is Almgren minimal. This is motivated by the attempt to classify all possible singularities for Almgren minimal sets--Plateau's problem in the setting of sets. In particular, a direct application of the above result leads to various types of new singularities for Almgren minimal sets, e.g.  the product of any paired calibrated cone (such as the cone over the $d-2$ skeleton of the unit cube in $\R^d, d\ge 4$) with homogeneous area minimizing hypercones (such as the Simons cone).

\medskip

\centerline {\large\textbf{R\'esum\'e.}} 

Dans cet article, on d\'emontre que le produit d'un ensemble calibr\'e appari\'e et d'un ensemble de codimension 1 calibr\'e par un calibration ''coflat'' avec un petit ensemble de singularit\'es est minimal au sens d'Almgren. Ceci est motiv\'e par la tentative de classifier toutes les singularit\'es possibles pour les ensembles minimaux d'Almgren--le probl\`eme de Plateau dans le cadre des ensembles. En particulier, une application directe du r\'esultat ci-dessus donne divers types de nouvelles singularit\'es pour les ensembles minimaux d'Almgren. Par exemple, le produit de tout c\^one calibr\'e appari\'e (comme le c\^one sur le squelette $d-2$ du cube unit\'e dans $\R^d, d\ge 4$) avec les hyperc\^ones homog\`ene minimisant la masse (comme le c\^one de Simons).

\medskip



\setcounter{section}{-1}

\section{Introduction}

In this paper, we prove that the product of a paired calibrated set and a calibrated set of codimension 1 calibrated by a coflat calibration with small singularity set is Almgren minimal. This is motivated by the attempt to classify all possible singularities for Almgren minimal sets--Plateau's problem in the setting of sets. In particular, a direct application of the above result leads to various types of new singularities for Almgren minimal sets.

Plateau's problem aims at understanding existence and local structure for physical objets that minimize the area while spanning a given boundary. A well known example is the soap films, which are objects of dimension 2 living in 3 dimensional ambient space. 

The approach of Plateau's problem involves the mathematical intepretation of the words ''objects, area, spanning''. Objects can be functionals (currents), measures (varifolds), $C^2$-manifolds (minimal surfaces), and sets (various notion of minimal sets), etc. They are different but closely related in many cases.

In case of soap films, the notion of (Almgren) minimal sets introduced by F.J. Almgren \cite{Al76} gives a very good descripton of the local behavior of soap films. In particular, the classification of singularities for 2-dimensional Almgren minimal sets in 3-dimensional ambient spaces (well known result of J. Taylor 1976 \cite{Ta}) coincides perfectly with what we can observe in soap film experiments. 

For general dimensions and codimensions, we know that any $d$-dimensional Almgren minimal set must be a manifold outside a $\H^d$-null set (See works of Almgren \cite{Al76}, David \& Semmes \cite{DS00}, David\cite{DEpi}). Points of this $\H^d$-null set are called singular points, and they do not admit tangent planes. A typical way to study local behavior around singular points is to look at the "tangent objects"-blow up limits at these points. It is proved by David \cite{DJT} that these blow-up limits are all minimal cones---minimal sets that are cones as well. Thus one tries to classify all possible minimal cones, or in other words, singularities.

In $\R^3$, the list of 2-dimensional minimal cones has been given by several mathematicians a century ago. (See for example \cite{La} or \cite{He}). They are, modulo isometry: a plane (which we also call a $\P$ set), a $\Y$ set (the union of 3 half planes that meet along a straight line where they make angles of 120 degree; this straight line is called its spine), and a $\T$ set (the cone over the 1-skeleton of a regular tetrahedron centered at the origin). See the pictures below. 

\centerline{\includegraphics[width=0.16\textwidth]{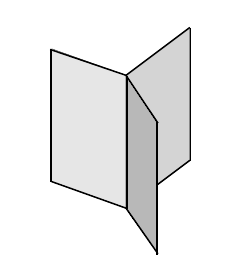} \hskip 2cm\includegraphics[width=0.2\textwidth]{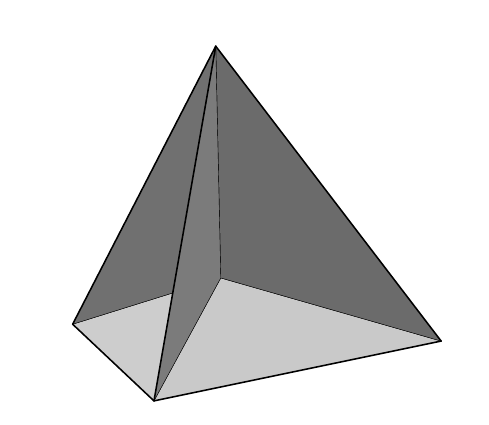}}
\nopagebreak[4]
\hskip 4cm a $\Y$ set\hskip 3.9cm  a $\T$ set\\
\nopagebreak[4]
\hskip 7cm Figure 1

In higher dimensions, even in dimension 4, the list of minimal cones is still very far from clear. 

Based on the above, one natural way to find new candidate for minimal cones, is by taking products of known minimal cones. In \cite{DJT}, David asked the question whether the product of 2 minimal cones stays minimal.

%
%

On the other hand, unlike currents, varifolds, etc., which are objects in dual spaces and inherit naturally algebraic structures such as dual, multiplicity, orientation, etc., minimal sets are just closed sets from their definition. Thus in general it is much more difficult to prove the Almgren minimality of a set. Furthermore, some properties that seem obvious in intuition are hard to prove for minimal sets (and hence minimal cones), such as
 the minimality for the product of two Almgren minimal sets. It is not even known whether the product of a general minimal set with a line is minimal. 
  
There are still partial results concerning the minimality for the product of two minimal sets:

(1) First, in \cite{topo} the authors proved that the product of $\R^n$ with any minimal set satisfying some additional property on the homology group of its complement is minimal. Since the product of a set with $\R^n$ is just a "thicker" version of the set, this product does not give essentially new topology. 

(2) In order to find new topological types of singularities, people turn to look at products of two minimal sets other than $\R^n$. This is apparently more complicated: a first attempt is the minimality of the products of two 1 or 2-dimensional $\Y$ sets (see Figure 1). Its minimality was conjectured by David in \cite{DEpi}, and proved by the author in \cite{YXY}. Compared to the previous result, this only gives the minimality of two particular minimal cones. But this gives a new minimal cone with completely new topology. The proof is very involved, relies heavily on the structure of $\Y$, and is not easy to generalise.

(3) The proof of the above result (2) is also based on the particular family of paired calibrations that calibrates $\Y$. The technique of paired calibrations has been introduced separately by Brakke \cite{Br91}, Lawlor and Morgan \cite{LM94} for proving the minimality of various sets of codimension 1. We then try to prove the minimality of the product of paired calibrated sets with other special minimal sets. The simplest candidate to take (other than $\R^n$) is probably the area-minimizing hypersurfaces. In \cite{product1}, the author proved the Almgren minimality for the product of a paired calibrated set and a big class of codimensional 1 area minimizing manifold (area minimizing manifold that admits smooth calibrations).

The above result in (3) only allows smooth manifolds as the second component of the product. Minimal cones other than $\R^n$ must have singularities, hence (3) does not give any new singularities (i.e. minimal cones) for Almgren minimal sets. 

In this paper, we take singular area minimizing sets into account. Let us first state one of our main results, which admits the result in (3) as a very special case:

\begin{thm} Let $n\ge 2,d\ge 2$ be integers. 

Let $U$ be an open bounded set in $\R^n$. Let $\b$ be an $n-2$-integral current without boundary, such that $B=spt(\b)\subset \pa U$.  Let $E\in C_1(\b,w,U)$, where $w$ is a coflat calibration in $U$, such that its singular set $S_w\subset E$, and $\H^{n-2}(S_w)=0$. 

Let $V\subset \R^d$ be a bounded strongly Lipschitz domain, so that $(\R^d,\pa V)$ admits a $C^1$ triangulation. Let $F\in PC(\mo,\nu)$, where $\mo=(\Omega_1,\cdots, \Omega_k)\in \cO_P(V,k)$ for some $k\in \N$, and $\nu\in \cV_c^k(\bar V)$ is such that $v_i-v_j$ does not vanish on $V\bs S_\nu$ for any $i<j$.

Suppose that $U\times V$ is a SR-domain. 

Then the product $E\times F$ is a $n+d-2$-dimensional Almgren minimal set in $U\times V$.
\end{thm}

Here $\cV_c^k(\bar V)$ stands for the class of coflat paired calibrations in a neighborhood of $\bar V$, and $S_\nu$ stands for the singular set of $v$. See the Definition \tb{3.3} for their precise definitions. Roughly speaking, the theorem says that if $E\subset \R^n$ is the support of a codimension 1 calibrated size-minimizing current so that its singular set is a $\H^{n-2}$ null subset of $E$, and $F$ is a paired calibrated set, then their product is Almgren minimal provided that their ambient domain admits some mild regularity.

The words "calibration" and "paired calibrations" in the theorem might be a little confusing because they look similar. So let us say a little bit more: 

Calibrated sets come from the theory of current. They are supports of a big class of size minimizing integral currents that span a given boundary, hence are minimal sets. The proof of their minimality is based on a dual argument: calibration forms. The theory of calibration dates back to \cite{HaLa82}, which applies first to manifolds, then to area minimizing currents.  See Section \tb{2} for detail.

On the other hand, the technique of paired calibrations has been introduced separately by Brakke \cite{Br91}, Lawlor and Morgan \cite{LM94} for proving the minimality of various sets of codimension 1. In contrast with the "spanning minimality" for calibrated sets, paired calibrated sets minimize the Hausdorff measure among classes of sets that satisfy some given separation condition. See Section \tb{3}.

It might be worth mentioning that many codimensional 1 calibrated sets are particular paired calibrated sets. See Example \tb{5.3} for a typical example. It follows that our result can also be applied to products of two codimensional 1 calibrated sets.

As pointed out before, singular minimal cones cannot be calibrated by ordinary smooth calibrations. In order to find new types of singular minimal cones, one has to consider calibrations with singularities. In \cite{HaLa82} the authors also introduced more general coflat calibrations--they are differential forms with singularities. Many area-minimizing hypercones admit coflat calibrations, see Example \tb{5.2}. As a result, by taking products, we obtain many new types of singularities for Almgren minimal sets. See Section \tb{5}.

Of course, the existence of singularities brings troubles. For example, we lose the smoothness--even the continuity--of the projection along the flow of the calibration, because the flow can touch singular sets. So the projection argument in \cite{product1}, or any other similar projection along calibration flow to calibrated sets does not work. Also, since we are in codimension higher than 1, we do not have natural orientation, separation, etc. at hand. This issue can be seen in the proof for the minimality of $\Y\times \Y$: it is already very involved even though $\Y$ is the simplest minimal cone with singularity. In this paper, we are dealing with general calibrated and paired calibrated sets with singularities, so more ideas are needed.

The general idea together with the organization of the article is the following: 

In section \tb{1} we introduce basic preliminaries and notations, and prove an approximation proposition which guarantees that in domains with some regularity (SR domains) we can restrict ourselves to recitifiable competitors in the proof of minimality.

In Section \tb{2 and 3} we discuss and prove the relations between the various minimality for codimensional 1 coflat calibrated sets and coflat paired calibrated sets. For calibrations with singularities, here we cannot do Lipschitz projection onto a calibrated set $E$ along the flow of the calibration form $w$. But we still manage to prove that the flow $\theta$ of $w$ that ends on the singular set of of $w$ or ends on the boundary of $E$ only touches $E$ at a $\H^{n-1}$-null set. As a result, given a competitor $F$ that satisfies some homological condition, we can project the part of $F$ in the regular region of $\theta$ along $\theta$, and the projection covers $E$ except for a $\H^{n-1}$-null set. Similar things work for coflat paired calibrations.

Section \tb{4} is the crucial part, in which we discuss various minimality results for the product of a codimension 1 calibrated set $E$ (with coflat calibration $w$) with singularities in $U\subset \R^n$ and a paired calibrated set $F$ (with a family of coflat calibrations $v_i,1\le i\le k$) in a domain $V\subset \R^d$. Given $A$ a competitor that satisfies the same homolgical conditions as $E\times F$, we would like to project $A$ to $E\times \bar V$ along the flow $\theta$ of $w$. But since the calibration $w$ has singularities, we can only project the part of $A$ inside the regular region of $\theta$. Of course this partial projection does not keep any homological property of $A$, so we have to "stitch" its image together in $E\times \bar V$ to get a set $A'\subset E\times \bar V$ that satisfy the desired homological property, with which we can still apply the paired calibrations $w\wg v_i,1\le i\le k$. This "stitch" has to be done carefully so that it costs us as small $\H^{n+d-2}$-measure as we want. Other operations are also needed, such as "flattening" the competitor $A$ down to $E\times \bar V$ near the boundary of $U\times V$ without costing too much measure. After all these operations, we can apply the paired calibrations $w\wg v_i,1\le i\le k$ to $A$, then pass by the intermediate set $A'$ we get the desired minimality of $A$.

In Section \tb{5} we give examples of families of codimensional 1 calibrated and paired calibrated cones to which we can apply our theorem. This gives several new families of singularities for Almgren minimal sets of codimension 2. A particular example is the product of the Simons cone with any paired calibrated set.

\textbf{Acknowledgement:} The author would like to thank Yongsheng Zhang for discussions and help on examples of coflat calibrated area minimizing hypercones. This work is supported by the National Natural Science Foundation of China (Grant No. 11871090, and Grant No. 12271018).

\section{Definitions and preliminaries}

\subsection{Basic notation and definitions}

$B(x,r)$ is the open ball with radius $r$ and centered on $x$;

$\overline B(x,r)$ is the closed ball with radius $r$ and center $x$;

In any metric space, $\H^d$ denotes the Hausdorff measure of dimension $d$;

A neighborhood of a set $E$ is just an open set that contains $E$;

For any set $E$ in a metric space, and any $r>0$, $B(E,r)=\{x:d(x,E)<r\}$, and $\bar B(E,r)=\{x:d(x,E)\le r\}$.

Let $Y\subset X$, then $i_{Y,X}:Y\to X$ denotes the inclusion map.

Let $C^S$, $C^\Delta$ denote the singular and simplicial chain complices. $H^S$, $H^\Delta$ denote singular homology and simplicial homology respectively.

\begin{defn}[reduced set] For every closed subset $E$ of $\R^n$, denote by
\be E_d^*=\{x\in E\ ;\ \H^d(E\cap B(x,r))>0\mbox{ for all }r>0\}\ee
 the $\H^d$-kernel of $E$. We say that $E$ is $d$ reduced if $E=E^*_d$.
\end{defn}

In the following context, for all $d\in \N$, all $d$-dimensional sets are supposed to be $d$-reduced.

\begin{defn}Let $d$ be an integer.
For any set $E$ and any point $x\in E$, let $\theta_d^*(E,x)$ and $\theta_{d,*}(E,x)$ denote the upper and lower $d$-density of $E$ at $x$ respectively:
$$\theta_d^*(E,x)=\limsup_{r\to 0}\frac{\H^d(E\cap B(x,r))}{\a_dr^d},\ \theta_{d,*}(E,x)=\liminf_{r\to 0}\frac{\H^d(E\cap B(x,r))}{\a_dr^d},$$
where $\a_d$ stands for the volume of the Euclidean unit $d$-ball. 

When $\theta_d^*(E,x)=\theta_{d,*}(E,x)$, we call it the $d$-density of $E$ at $x$ and denote it by $\theta_d(E,x)$.
\end{defn}

If $E$ is a $d$-rectifiable set in $\R^n$, denote by $T_xE\in G(n,d)$ the approximate tangent plane (if it exists and is unique) of $E$ at $x$.

Note that when $E$ is rectifiable, then for $\H^d$-a.e. $x\in E$, $T_xE$ and $\theta(E,x)$ exist, and $\theta(E,x)=1$. 

For any $d$-integral current $T$ in $\R^n$, let $\spt(T)$ denote its support. It is $d$-reduced.

\begin{defn}[cf. \cite{DP07} Section 3.3] Let $A\subset \R^n$. For any non negative integer $d$, let $C^I_d(A)=C^I_d(A,\Z)$ denote the class of all $d$-integral currents on $A$. Let $\pa_d: C^I_d(A)\to C^I_{d-1}(A)$ be the boundary operator. The $d$-dimensional integral rectifiable homology group $H^I_d(A)$ is the quotient 
\be H^I_d(A)=\frac{Ker \pa_d}{Im \pa_{d+1}}.\ee
\end{defn}

\begin{defn}Let $d$ be an integer. Let $X$ be a topological space $X$, $G$ be an abelien group. Let $H$ be a homology theory ($H=H^S, H^\Delta or H^I$) that can be defined on $X$ with coefficients in $G$. 

$1^\circ$ Let $\Gamma$ be an element in the associated chain complex $C_d(X,G)$, then $\lg \Gamma\rg_{H_d(X,G)}$ denotes the class of all chains in $C_d(X,G)$ that are homologic to $\Gamma$: $ \lg \Gamma\rg_{H_d(X,G)}=\{\Gamma'\in C_d(X,G):$ there  exists $R\in C_{d+1}(X,G)$ so that  $\pa R=\Gamma-\Gamma'\}$.

$2^\circ$ We write $\Gamma\sim_{H_d(X,G)}\Gamma'$ if the two chains are homologic.

$3^\circ$ Let $f:X\to Y$ be continuous (Lipschitz when $C=C^I$ and $H=H^I$). Let $C_{d,G}(f): C_d(X,G)\to C_d(Y,G)$ and $H_{d,G}(f):H_d(X,G)\to H_d(Y,G)$ be the homomorphisms induced by $f$. When $f$ is Lipschitz, we also let $f_\sharp$ denote $C^I(d,\Z)(f)$ for short.
\end{defn}

Given a subset $X$ of a an Euclidean space, and a rectifiable set $\Gamma\subset X$ with a prescribed orientation, $[\Gamma]_{C^I_d}$ denotes the induced integral current; Given a $d$-simplicial $G$-chain $\Gamma$, $[\Gamma]_{C^I_d}$ (when $G=\Z$) and $[\Gamma]_{C^S_{d,G}}$ denote the induced integral current and singular chain respectively; and when $\pa \Gamma=0$, let $[\Gamma]_{H^\Delta_d(X,G)}$ and $[\Gamma]_{H^S_d(X,G)}$ denote the induced element in $H^\Delta_d(X,G)$ and $H^S_d(X,G)$ respectively. If $\Gamma\subset X$ is an oriented $C^1$ manifold, and $(X,\Gamma)$ is triangulable, then different triangulations on $(X,\Gamma)$ give different but homological elements in $C^\Delta_d(X,G)$, which induce a same element $[\Gamma]_{H^\Delta_d(X,G)}$ in the simplicial homology group and $[\Gamma]_{H^S_d(X,G)}$ in the singular homology group provided $\Gamma$ has no boundary.

%
%
%
%
%

\begin{defn}Let $d\le n$. Let $A,B$ be subsets of $\R^n$. We say that $A$ and $B$ are $d$-essentially disjoint, if $\H^d(A\cap B)=0$.
\end{defn}

\begin{defn}Given a $d$-rectifiable set $E\subset \R^n$, and a measurable unit tangent $d$-vector field $w$ on $E$, let $\H^d\lf_E\wg w$ be the induced rectifiable current: for any smooth $d$-form $\varphi$ with compact support on $\R^n$,
\be \lg\H^d\lf_E\wg w, \varphi\rg=\int_E\lg w(x), \varphi(x)\rg d\H^d(x).\ee
\end{defn}

\begin{defn}Let $d<n$ be integers. Let $E\subset \R^n$ be a closed set. We say that 

$1^\circ$ $E$ is $d$-integral regular, if there exists a $d$-integral current $T$ so that $E_d^*=\spt(T)$;

$2^\circ$ $E$ is $d$-simplicial regular, if there exists a $C^1$ simplicial complex $\cK$, so that $E$ is a finite union of faces of dimension no more than $d$ in $\cK$.

For any class of sets $\F$, set 
\be\F^{IR,d}=\{F\in \F: F\mbox{ is }d\mbox{-integral regular}\},\ee
\be\F^{SR,d}=\{F\in \F: F\mbox{ is }d\mbox{-simplicial regular}\},\ee
and 
\be\F^{R,d}=\{F\in \F: F\mbox{ is }d\mbox{-rectifiable with locally finite }\H^d\mbox{ measure}\}.\ee
\end{defn}

Clearly we have, for any $d$, 
\be \F^{SR,d}\subset\F^{IR,d}\subset \F^{R,d}.\ee

\subsection{Minimal set, mass minimizing current, relation}

\begin{defn}[General definition of minimal sets]Let $0<d<n$ be integers. Let $U\subset \R^n$ be an open set. A relatively closed set $E\subset U$ is said to be minimal of dimension $d$ in $U$ with respect to the competitor class $\mathscr F$ (which contains $E$) if 
\be \H^d(E\cap B)<\infty\mbox{ for every compact ball }B\subset U,\ee
and
\be \H^d(E\bs F)\le \H^d(F\bs E)\ee
for any competitor $F\in\mathscr F$.
\end{defn}

\begin{defn}[Almgren competitor (Al competitor for short)] Let $U$ be an open subset of $\R^n$. Let $E$ be a closed subset in $\bar U$. An Almgren competitor for $E$ in $U$ is a closed set $F\subset \bar U$ that can be written as $F=\varphi_1(E)$, where $\varphi_t:\bar U\to \bar U,t\in [0,1]$ is a family of continuous mappings such that 
\be \varphi_0(x)=x\mbox{ for }x\in \bar U;\ee
\be\mbox{ the mapping }(t,x)\to\varphi_t(x)\mbox{ of }[0,1]\times \bar U\mbox{ to }\bar U\mbox{ is continuous;}\ee
\be\varphi_1\mbox{ is Lipschitz,}\ee
  and if we set $W_t=\{x\in U\ ;\ \varphi_t(x)\ne x\}$ and $\widehat W=\bigcup_{t\in[0.1]}[W_t\cup\varphi_t(W_t)]$,
then
\be \widehat W\mbox{ is relatively compact in }U.\ee
 
Such a $\varphi_1$ is called a deformation in $U$, and $F$ is also called a deformation of $E$ in $U$.

The class of all Almgren competitors of $E$ in $U$ is denoted by $\F_{Al}(E,U)$.
\end{defn}

Note that by continuity of $\varphi_t$, we know that $\varphi_t|_{\pa U}=id, t\in [0,1]$.

\begin{defn}[Almgren minimal sets]
Let $0<d<n$ be integers, $U$ be an open set of $\R^n$. A $d$-dimensional Almgren-minimal set $E$ in $U$ is a minimal set defined in Definition \tb{1.8} while taking the competitor class $\mathscr F$ to be the class of all Almgren competitors for $E$.\end{defn}

In this paper, we mainly consider Almgren minimal sets. So in the following, minimal set means Almgren minimal set unless otherwise specified.

\begin{defn}[Integral homology competitor]Let $n,d$ be integers, $d<n$. Let $U\subset \R^n$ be an open set. Let $E\subset\bar U$ be a closed set. Set $B=E\cap \pa U$. Let $H$ be a class of $d$-integral currents supported in $E$. A set $A\subset \bar U$ is called a $(d,B,H,U)$-integral homology competitor of $E$, if $A\cap \pa U=B$, and for each element $S\in H$, there exists a $d$-integral current $T$ supported in $A$, so that $T\sim_{H^I_d(U\cup (spt(S)\cap \pa U), \Z)}S$. The class of all $(d,B,H,U)$-integral homology competitors of $E$ is denoted by $\F_{ihc}(E,d,B,H,U)$. \end{defn}

\begin{defn}[Integral current spanner] Let $n,d$ be integers, $d<n$. Let $U\subset \R^n$ be an open set. Let $B\subset \pa U$ be a closed set. Let $H$ be a set of $d-1$-integral currents supported in $B$ without boundary. We say that a closed set $E\subset \bar U$ is a $(d,B,H,U)$-integral current spanner, if $E\cap \pa U=B$, and $(i_{B,E})_\sharp(H)\subset \lg 0\rg_{C^I_{d-1}(E)}$, That is, for any $S\in H$, there exists a $d$-integral current $T$ supported in $E$, so that $\pa T=S$. The class of all $(d,B,H,U)$-integral current spanners is denoted by $\F_{ics}(d,B,H,U)$.
\end{defn}

\begin{rem}Let $d,n, U, E,B, H$ be as in the definition of integral homology competitor. Suppose $H$ satisfies that : for each $S\in H$, $\pa S$ is supported in $B$. Let $H'=\{\pa S:S\in H\}$. Then it is easy to see that
\be \F_{ihc}(E,d,B,H,U)\subset \F_{ics}(d,B,H',U).\ee
\end{rem}

\begin{defn}[Spanning competitors] Let $n,d$ be integers, $d<n$. Let $G$ be an abelian group. Let $U\subset \R^n$ be an open subset. Let $B\subset \pa U$ be a closed set. Let $H$ be a subset of $H^S_{d-1}(B,G)$. We say that a closed set $E\subset \bar U$ is a $(d,G,B,H,U)$-spanning competitor, or a $(d,G,B,H,U)$-spanner, if $E\cap \pa U= B$, and the map $H^S_{d-1,G}(i_{B,E}): H^S_{d-1}(B,G)\to H^S_{d-1}(E,G)$ satisfies that $H\subset Ker H^S_{d-1,G}(i_{B,E})$. When the dimension $d$ is already fixed, we also say that $E$ spans $H$ in $U$ with coefficient in $G$. Moreover, if $G=\Z$, we also say that $E$ spans $H$ in $U$.

Denote by $\F_{sc}(d,G,B,H,U)$ the class of all $(d,G,B,H,U)$-spanners. 
\end{defn}

\begin{defn}[Topological competitors]Let $n,d$ be integers, $d<n$. Let $G$ be an abelian group. Let $U\subset \R^n$ be an open set. Let $B\subset \pa U$ be closed. Let $H$ be a subset of $H^S_{n-d-1}(\pa U\bs B, G)$. We say that a closed set $E\subset \bar U$ is a $(d,G,B,H,U)$-topological competitor, if $E\cap \pa U=B$, and the map $H^S_{n-d-1,G}(i_{\pa U\bs B,\bar U\bs E}):H^S_{n-d-1}(\pa U\bs B, G)\to H^S_{n-d-1}(\bar U\bs E, G)$ satisfies that $H\cap ker H^S_{n-d-1,G}(i_{\pa U\bs B,\bar U\bs E})\bs\{0\}=\emptyset$. 

Denote by $\F_{tc}(d,G,B,H,U)$ the class of all $(d,G,B,H,U)$-topological competitors.
\end{defn}

The following proposition tells us that all the above classes are stable under deformations in $U$.

\begin{pro}Let $n,d$ be integers, $d<n$. Let $G$ be an abelian group. Let $U$ be an open subset of $\R^n$. Let $B\subset \pa U$ be closed. Let $E$ be closed in $\bar U$ so that $E\cap \pa U=B$. Let $H_1$ be a class of $d$-integral currents supported in $E$. Let $H_2$ be a set of $d-1$-integral currents supported in $B$ without boundary. Let $H_3$ be a subset of $H^S_{d-1}(B,G)$.  Let $H_4$ be a subset of $H^S_{n-d-1}(\pa U\bs B, G)$. Then for $\F=\F_{ihc}(E,d,B,H_1,U), \F_{ics}(d,B,H_2,U), \F_{sc}(d,G,B,H_3,U)$, and for $\F_{tc}(d,G,B,H_4,U)$ when $(\R^n,\pa U)$ admits a $C^1$ triangulation,
\be E\in \F\Rightarrow \F_{Al}(E,U)\subset \F.\ee
In particular, if $E$ is minimal of dimension $d$ with respect to $\F$, then $E$ is Almgren minimal of dimension $d$ in $U$.
\end{pro}

\nd The proof of the cases for $\F=\F_{ics}(d,B,H_2,U), \F_{sc}(d,G,B,H_3,U)$ and $\F_{tc}(d,G,B,H_4,U)$ can be found directly in Proposition 2.14 of \cite{product1}. A small remark is that in the statement of the cited Proposition 2.14 of \cite{product1}, we asked $\pa U$ to be a $C^1$ manifold to guarantee that $\F_{tc}$ contains all Almgren competitors. Here we only ask that $(\R^n,\pa U)$ admits a $C^1$ triangulation, because one can see directly from the proof of Proposition 2.14 of \cite{product1}, that this is the only property we used. 

Let us verify the case for $\F=\F_{ihc}(E,d,B,H_1,U)$: let $A\in \F_{Al}(E,U)$. By definition there exists a Lipschitz deformation $\varphi_t:\bar U\to \bar U$ as in Definition \tb{1.9}, so that $A=\varphi_1(E)$.

Now for each $S\in H_1$, let $T=\varphi_{1,\s}(S)$, and $R=\varphi_\s([0,1]\times S)$ where $[0,1]$ is endowed with the positive orientation. Then they are both supported in $U\cup (spt(S)\cap \pa U)$, and $\pa R=S-T$. That is, $T$ is a $d$-integral current supported in $A$, so that $T\sim_{H^I_d(U\cup (spt(S)\cap \pa U), \Z)}S$. Hence $A\in \F=\F_{ihc}(E,d,B,H_1,U)$.
\qed

In the following we wish to give some approximation theorem, which says that regular sets are ''dense'' in measure in the above classes. For this, we need some regularity of the domain. Briefly, we ask that the domain can "retract" inside itself gradually in a Lipschitz way.

\begin{defn}[$(L,\e)$-self retract domain] Let $L\ge 1$, $\e>0$. Let $U$ be an open subset of $\R^n$. A map $\varphi: \bar U\times [0,1]\to \bar U$ is called a $(L,\e)$-self retract of $U$, if it satisfies:

(1) $\varphi$ is $L$-Lipschitz;

(2) If we denote $\varphi_t=\varphi(\cdot, t)$, then $\varphi_0=id$;

(3) For every $t\in (0,1]$, $\varphi_t(\bar U)\subset U$, and $\varphi_t|_{U^-_{\e}}=id$. Here for each $\d>0, U^-_\d:=U\bs B(\pa U, \d)$;

(4) $\varphi|_{\pa U\times [0,1]}$ is a $C^1$-diffeomorphism.

For $1\le k\le\infty$, we say that $\varphi$ is a $C^k$ $(L,\e)$-self retract of $U$, if in addition $\varphi$ is of class $C^k$, and $\varphi|_{\pa U\times [0,1]}$ is a $C^k$ diffeomorphism.

We say that an open subset $U$ of $\R^n$ is a $(L,\e)$-self retract domain, $(L,\e)$-SR domain for short, if for each $\d\le \e$, there exists a $(L,\d)$-self retract of $U$. It is called a $C^k$ $(L,\e)$-SR domain, if it admits a $C^k$ $(L,\d)$-self retract for each $\d\le \e$. We say that $U$ is an SR-domain (resp. $C^k$ SR-domain), if it is a $(L,\e)$-SR domain (resp. $C^k$ $(L,\e)$-SR domain) for some $L\ge 1$ and $\e>0$.
\end{defn}

\begin{rem}$1^\circ$ The above definition for SR-domain is relatively weak, because we only need it to guarantee the approximation theorem below. Some other natural properties that are also suitable for the name "retract" are not imposed for SR-domains. For example, we do not ask that $\varphi_t$ is a homeomorphism for each $t$, and we do not ask that $\varphi_t(\bar U)\subset \varphi_s(\bar U)$ provided $t\ge s$, or $\varphi_1(\bar U)=U^-_\e$, etc. 

With such a weak condition, it is not clear if the product of two SR domains is always an SR domain. But it is, if we add some additional property, e.g. if each $\varphi_t$ is a $C^1$ diffeomorphism. Also, the strongly Lipschitz domains (used in \cite{product1} to get the same approximation theorem) are SR domains, and the product of two strongly Lipschitz domains are still strongly Lipschitz. 

An important class of strongly Lipschitz domains is the class of convex domains: in the study of minimal cones, one often take convex domains, e.g., convex hull of a minimal cone, or the unit ball. 
\end{rem}

\begin{pro}[Approximation]Let $n,d$ be integers, $d<n$. Let $G$ be an abelian group. Let $U$ be a bounded SR domain in $\R^n$. Let $B\subset \pa U$ be closed so that $\H^{d-1}(B)<\infty$. Let $H_1$ be a class of $d$-integral currents supported in $E$. Let $H_2$ be a set of $d-1$-integral currents supported in $B$ without boundary. Let $H_3$ be a subset of $H^S_{d-1}(B,G)$.  Let $H_4$ be a subset of $H^S_{n-d-1}(\pa U\bs B, G)$. Then

$1^\circ$ If $B$ is $d-1$-rectifiable, then for $\F=\F_{ihc}(E,d,B,H_1,U), \F_{ics}(d,B,H_2,U), \F_{sc}(d,G,B,H_3,U)$, and for $\F_{tc}(d,G,B,H_4,U)$ when $(\R^n,\pa U)$ admits a $C^1$ triangulation, we have
\be \inf\{\H^d(F),F\in \F\}=\inf\{\H^d(F),F\in \F^{R,d}\};\ee

$2^\circ$ If $B$ is $d-1$-integral regular, then for $\F=\F_{ihc}(E,d,B,H_1,U), \F_{ics}(d,B,H_2,U), \F_{sc}(d,G,B,H_3,U)$, and for $\F_{tc}(d,G,B,H_4,U)$ when $(\R^n,\pa U)$ admits a $C^1$ triangulation, we have
\be \inf\{\H^d(F),F\in \F\}=\inf\{\H^d(F),F\in \F^{IR,d}\};\ee

$3^\circ$ If $B$ is $d-1$-simplicial regular, then for $\F=\F_{ihc}(E,d,B,H_1,U)$, $\F_{ics}(d,B,H_2,U)$, $\F_{sc}(d,G,B,H_3,U)$, and for $\F_{tc}(d,G,B,H_4,U)$ when $(\R^n,\pa U)$ admits a $C^1$ triangulation, we have
\be \inf\{\H^d(F),F\in \F\}=\inf\{\H^d(F),F\in \F^{SR,d}\}.\ee
\end{pro}

\nd First, for $\F=\F_{ihc}(E,d,B,H_1,U), \F_{ics}(d,B,H_2,U), \F_{sc}(d,G,B,H_3,U)$, and for $\F_{tc}(d,G,B,H_4,U)$, if $\inf\{\H^d(F),F\in \F\}=\infty$, then there is nothing to prove. Hence we suppose that $\inf\{\H^d(F),F\in \F\}<\infty$.

To prove Proposition \tb{1.19}, for each given closed set $A\subset \bar U$ with finite $\H^d$ measure such that $A\cap \pa U=B$, and any given $\e>0$, we will construct a set $F$ so that $\H^d(F)<\H^d(A)+\e$, $F$ is in the same competitor class as $A$, and $F$ is as regular as $B$.

Since $U$ is a SR domain, by definition, there exists $L\ge 1, \e_0>0$, so that for each $\d\le \e_0$, there exists a $(L,\d)$-self retract of $U$.

Since $\H^{d-1}(B)<\infty$, we know that $0=\H^d(B)=\lim_{r\to 0}\H^d(A\cap B(\pa U, r))$, we can take $\d<\e_0$ so that $\H^d(A\cap B(\pa U, \d))<\e/3L^d$.

Let $\varphi$ be a $(L, \d)$-self retract of $U$. Let $t_0=\frac{\e}{3L^d\H^{d-1}(B)}$. Let $\psi=\varphi_{t_0}$, and set  $F_0=\psi(A)\cup \varphi(B\times [0,t_0])$. Then we have
\be \begin{split}\H^d(F_0)&\le \H^d(\psi(A\cap U^-_\d))+\H^d(\psi(A\cap B(\pa U,\d)))+\H^d(\varphi(B\times[0,t_0]))\\
&\le \H^d(A\cap U^-_\d)+L^d\H^d(A\cap B(\pa U,\d))+L^d\H^d(B\times[0,t_0])\\
&\le \H^d(A)+\frac\e3+\frac\e3= \H^d(A)+\frac{2\e}{3}.
\end{split}\ee

Now we will deform the set $F_0$ in a polyhedral complex to get it regular. Let $d_0=d(\psi(\bar U), \pa U)$, let $A_0=\psi(\bar U)$, let $A_1=\bar B(A_0, \frac{d_0}{2})\subset U_{\frac{d_0}{2}}^-$.

Next we apply Lemma 5.2.6 of \cite{Fv}, to the set $F_0$, the domain $U$, the compact set $A_1$, and get a deformation $g$ in $U$, and a polyhedral complex $\cK$ of diameter less than $10^{-1}d_0$, so that 

$1^\circ$ $A_1\subset |\cK|^\circ\subset |\cK|\subset U$;

$2^\circ$ $g|_{U\bs |\cK|}=id$, $g(|\cK|)\subset |\cK|$;

$3^\circ$ $g(F_0)\cap |\cK|^\circ$ is a union of faces of dimension no more than $d$ of $\cK$;

$4^\circ$ $\H^d(g(F_0))\le \H^d(F_0)+\e/3$.

Moreover, after the proof of Lemma 5.2.6 of \cite{Fv}, we know that for each $n-1$-face $\sigma$ of $\pa|\cK|$, if $\Sigma$ is the polyhedron of $\cK$ so that $\sigma$ is a face of $\Sigma$, then there exists a point $\xi\in \Sigma\bs F_0$, so that $g(F_0)\cap \sigma=\pi_\xi(F_0\cap \Sigma)\cap \sigma$, where $\pi_\xi: \Sigma\bs \{\xi\}\to \pa\Sigma$ is the central projection. Note that since the diameters of the polygons in $\cK$ are less than $10^{-1}d_0$, and $A_0\subset A_1\subset |\cK|^\circ$, we know that $\Sigma\cap A_0=\emptyset$. Hence $F_0\cap \Sigma=\varphi(B\times [0,t_0])\cap \Sigma$. Hence \be g(F_0)\cap \sigma=\pi_\xi(\varphi(B\times [0,t_0]))\cap \Sigma).\ee

As a result, since $\varphi|_{\pa U\times [0,1]}$ is a $C^1$-diffeomorphism, if $B$ is $d-1$-rectifiable (resp. integral regular, simplicial regular ), then $g(F_0)\cap \sigma$ is $d$-rectifiable (resp. integral regular, simplicial regular) for any $n-1$-face $\sigma$ of $\pa|\cK|$. Hence $g(F_0)\cap \pa|\cK|$ is also $d$-rectifiable (resp. integral regular, simplicial regular).

Recall that $g(F_0)\bs |\cK|=F_0\bs |\cK|\subset F_0\bs A^\circ_1\subset \varphi(B\times [0,t_0])$, which is $d$-rectifiable (resp. integral regular, simplicial regular) as long as $B$ is $d-1$-rectifiable (resp. integral regular, simplicial regular). Also, by $3^\circ$, the part $g(F_0)\cap |\cK|^\circ$ is $d$-simplicial regular. 

Set $F=g(F_0)$. Then the union $F=g(F_0)=[g(F_0)\cap |\cK|^\circ]\cup[g(F_0)\cap \pa|\cK|]\cup [g(F_0)\bs |\cK|]$ satisfies that
\be \begin{split}F\mbox{ is }d\mbox{-rectifiable (resp. integral regular, simplicial regular) as long as }
B\mbox{ is }d-1\\\mbox{-rectifiable (resp. integral regular, simplicial regular)}.\end{split}\ee

Moreover, by \tb{$4^\circ$} and \tb{(1.19)} we have
\be \H^d(F)\le \H^d(F_0)+\frac\e3\le \H^d(A)+\frac{2\e}{3}+\frac\e3=\H^d(A)+\e.\ee

The rest is to prove that for $\F=\F_{ihc}(E,d,B,H_1,U)$, $\F_{ics}(d,B,H_2,U)$, $\F_{sc}(d,G,B,H_3,U)$, and for $\F_{tc}(d,G,B,H_4,U)$ when $(\R^n,\pa U)$ admits a $C^1$ triangulation, $A\in \F\Rightarrow F\in \F$. Note that by Proposition \tb{1.16}, the classes $\F=\F_{ihc}(E,d,B,H_1,U)$, $\F_{ics}(d,B,H_2,U)$, $\F_{sc}(d,G,B,H_3,U)$, and the class $\F=\F_{tc}(d,G,B,H_4,U)$ when $(\R^n,\pa U)$ admits a $C^1$ triangulation, are all stable under deformations in $U$, and $F=g(F_0)$ where $g$ is a deformation in $U$, hence it is enough to prove that $F_0$ is in $\F$ provided $A$ is in $\F$.

--$\F=\F_{ihc}(E,d,B,H_1,U)$: Let $A\in \F$. By definition for each $S\in H_1$,  there exists a $d$-integral current $T$ supported in $A$, so that $T\sim_{H^I_d(U\cup (spt(S)\cap \pa U), \Z)}S$. Let $R=(-1)^d\varphi_\s(T\times [0,t_0])$ where $[0,t_0]$ is endowed with the positive orientation. Then we have $spt(R)\subset U\cup (spt(S)\cap \pa U)$, and 
\be\pa R=\varphi_\s([T\times \{t_0\}-T\times \{0\}]+(-1)^d\pa T\times [0,t_0])=\varphi_\s(T\times \{t_0\}+(-1)^d\pa T\times [0,t_0])-T.\ee
Let $T'=\varphi_\s(T\times \{t_0\}+(-1)^d\pa T\times [0,t_0])$, then it is supported in $F_0$ by definition, and 
\be T'\sim_{H^I_d(U\cup (spt(S)\cap \pa U), \Z)} T\sim_{H^I_d(U\cup (spt(S)\cap \pa U), \Z)}S.\ee
As a result, we know that $F_0\in \F$.


-- $\F=\F_{ics}(d,B,H_2,U)$: Let $A\in \F$. By definition for each $S\in H_2$, there exists a $d$-integral current $T$ supported in $A$ so that $\pa T=S$. Then $\psi_\sharp (T)+\varphi_\sharp (S\times [0,t_0])$ is a $d$-integral current supported in $F_0$ with boundary $S$. Hence $F_0\in \F$.

-- $\F=\F_{sc}(d,G,B,H_3,U)$: Let $A\in \F$. We would like to prove that $F_0\in \F$. So let $\theta$ be a $d-1$-$G$-singular cycle in $B$ which represents an element in $H_3$. Since $A\in \F$, we know that there exists a $d$-$G$-singular cycle $\Theta$ in $A$ so that $\pa \Theta=\theta$.

Now for the map $\psi$ that maps $A$ to $\psi(A)$, we know that $\pa C^S_{d,G}(\psi)(\Theta)= C^S_{d-1,G}(\psi)(\theta)$. Since $ C^S_{d,G}(\psi)(\Theta)$ is supported in $F_0$, $C^S_{d-1,G}(\psi)(\theta)$ represents a zero element in $H^S_{d-1}(F_0,G)$.

On the other hand, we have $\varphi:B\times [0,t_0]\to \varphi(B\times [0,t_0])$ is a diffeomorphism, and $\varphi(B\times [0,t_0])\subset F_0$, hence $C^S_{d-1,G}(\varphi(\cdot, 0))(\theta)\sim_{H^S_{d-1}(F_0,G)}C^S_{d-1,G}(\varphi(\cdot,t_0))(\theta)$ in $F_0$. Note that $\varphi(\cdot, 0)|_B=id$ and $\varphi(\cdot, t_0)|_B=\psi|_B$, hence $0\sim_{H^S_{d-1}(F_0,G)}C^S_{d-1,G}(\psi)(\theta)\sim_{H^S_{d-1}(F_0,G)}\theta=C^S_{d-1,G}(i_{B,F_0})(\theta)$ in $F_0$, which yields that $H^S_{d-1,G}(i_{B,F_0})(\lg \theta\rg_{H^S_{d-1}(F_0,G)})=0$ in $H_{d-1}(F_0,G)$. 

-- $\F=\F_{tc}(d,G,B,H_4,U)$: Let $A\in \F$. We would like to prove that $F_0\in \F$. We prove by contradiction, so suppose that $F_0\not\in\F$. By definition, $H_4\cap ker H^S_{n-d-1,G}(i_{\pa U\bs B,\bar U\bs F_0})\bs\{0\}\ne\emptyset$. Since $(\bar U,\pa U)$ admits a $C^1$ triangulation, and $F_0$ is closed, we know that $\bar U\bs F_0$ admits a locally finite $C^1$-triangulation $\cT_1$, so that there exists a simplicial $n-d$-$G$-chain $\Gamma$ of $\cT_1$, so that $\gamma=\pa\Gamma$ is in $\pa U\bs F_0$, and $\gamma$ respresents an element in $H_4\cap ker H^S_{n-d-1,G}(i_{\pa U\bs B,\bar U\bs F_0})\bs\{0\}$.

Let $\d_0=d(spt(\Gamma), F_0)$. 

Now we know that $\varphi:\pa U\times [0,t_0]$ is a $C^1$-diffeomorphism, and $(\R^n, \pa U)$ admits a $C^1$-triangulation, we know that $M:=\varphi(\pa U\times \{t_0\})$ is $n-1$-simplicial regular. By the transversality theorem, we can suppose that $\Gamma$ is transverse to $M$. Again since the map $\varphi_{t_0}:\bar U\to U$ satisfies that $\varphi_{t_0}|_{\pa U}$ is transverse to $\Gamma$, we can find a map $g:\bar U\to U$,  so that $g|_{\pa U}=\varphi_{t_0}|_{\pa U}$, $g$ is transverse to $\Gamma$, and $||g-\varphi_{t_0}||_\infty<10^{-1}\d_0$.

Let $\Gamma_1=\Gamma|_{\varphi(\pa U\times [0,t_0])}$, and let $\Gamma_2=\Gamma-\Gamma_1$. Since $\varphi(\pa U\times [0,t_0])$ is simplicial regular, whose boundary is $\varphi(\pa U\times \{t_0\})\cup \varphi(\pa U\times \{0\})=M\cup \pa U$, both are transverse to $\varphi_{t_0}$, hence $\Gamma_1$ is a well defined simplicial $n-d$-$G$-chain. Let $\Theta_1=\varphi^{-1}(\Gamma_1)$ (see \cite{topo} Proposition 2.36 for the definition and the boundary property of the inverse image of a simplicial chain under a transversal map). Let $\Theta_3=C^\Delta_{n-d, G}(\pi)(\Theta_1)$, where $\pi: [0,t_0]\times \pa U\to \pa U$ is the projection to the second component.

By definition of $\Gamma_1$ and $\Gamma_2$, we know that $spt(\Gamma_2)\subset \bar U\bs \varphi(\pa U\times [0,t_0])\subset g(\bar U)$. Since $g$ is transverse to $\Gamma$, $\Theta_2:=g^{-1}(\Gamma_2)$ is a $n-d$-$G$-chain in $\bar U$.

Let $\Theta=\Theta_3+\Theta_2$.

Let us first prove that $\pa\Theta=\gamma$. In fact, let $\gamma_1=\pa \Gamma_1$, then it is supported in $\pa M\cup \pa U$. Let $\gamma_2=\gamma-\gamma_1$. Then it is supported in $M$, and since $\Gamma=\Gamma_1+\Gamma_2$, we know that $\pa \Gamma_2=\gamma_2$.

Then we have 
\be\begin{split} \pa \Theta_3&=\pa C^\Delta_{n-d, G}(\pi)\circ\varphi^{-1}(\Gamma_1)=C^\Delta_{n-d, G}(\pi)\circ\varphi^{-1}(\pa\Gamma_1)=C^\Delta_{n-d, G}(\pi)\circ\varphi^{-1}(\gamma_1)\\
&=C^\Delta_{n-d, G}(\pi)\circ\varphi^{-1}(\gamma-\gamma_2)
=C^\Delta_{n-d, G}(\pi)\circ\varphi^{-1}(\gamma)-C^\Delta_{n-d, G}(\pi)\circ\varphi^{-1}(\gamma_2)\\
&=\gamma-\varphi_{t_0}^{-1}(\gamma_2).
\end{split}\ee

As for $\Theta_2$, we have
\be \pa \Theta_2=\pa g^{-1}(\Gamma_2)=g^{-1}(\pa\Gamma_2)=\varphi_{t_0}^{-1}(\gamma_2)\ee
since $g|_{\pa U}=\varphi_{t_0}|_{\pa U}$ and $\gamma_2$ is supported in $M=\varphi_{t_0}(\pa U)$.

As a result, we have
\be \pa \Theta=\pa \Theta_3+\pa\Theta_2=\gamma.\ee

Let us verify that $spt(\Theta)\cap A=\emptyset$. 

In fact, for any $x\in A$, since $\varphi_{t_0}(A)\subset F_0$, 
\be \begin{split}d(g(x), spt(\Gamma))&\ge d(\varphi_{t_0}(x), spt(\Gamma))-d(g(x), \varphi_{t_0}(x))\\
&\ge d(F_0, spt(\Gamma))-10^{-1}\d_0\ge \frac12\d_0,\end{split}\ee
and thus $g(x)\not\in \Gamma$, which implies $x\not\in spt(\Theta_2)$, hence $spt(\Theta_2)\cap A=\emptyset$; 

On the other hand, for any $(x,t)\in B\times [0,t_0]$, we know that 
\be d(\varphi(x,t), spt(\Gamma))\ge d(F_0, spt(\Gamma))\ge \d_0,\ee
thus $\varphi(B\times [0,t_0])\cap spt(\Gamma_1)=\emptyset$, which implies that $(B\times [0,t_0])\cap \varphi^{-1}(spt(\Gamma_1))=\emptyset$. But by definition, $spt(\Theta_1)\subset \varphi^{-1}(spt(\Gamma_1))$, hence $(B\times [0,t_0])\cap spt(\Theta_1)=\emptyset$, and therefore $B\cap \pi(spt(\Theta_1))=\emptyset$. Again by definition, $spt(\Theta_3)\subset \pi(spt(\Theta_1))$, we get $spt(\Theta_3)\cap B=\emptyset$. Since $spt(\Theta_3)\subset \pa U$, and $B=A\cap \pa U$, we have $spt(\Theta_3)\cap A=\emptyset$.

Altogether we get $spt(\Theta)\cap A=\emptyset$. But this contradicts the fact that $A\in \F$.\qed

\section{Minimality for 1 codimensional calibrated sets of multiplicity 1}

We fix an orthonormal basis $\{e_1,\cdots, e_n\}$ of $\R^n$. Let $\{e_1^*,\cdots, e_n^*\}$ denote the dual basis in $\bigwedge^1(\R^n)$. For any $d\le n$, and any $d$-vector $\xi\in \wg_d(\R^n)$, let $\xi^*\in \wg^d(\R^n)$ denote its dual $d$-covector. Also, for any $d$-covector $\zeta\in \wg^d(\R^n)$, let $\zeta_*\in \wg_d(\R^n)$ denote its dual vector.

\begin{defn}$1^\circ$ Let $W$ be an open subset of $\R^n$. A (measurable) $n-1$-form $w$ in $W$ is called a coflat $n-1$ form in $W$, if there exists a relatively closed subset $S_w$ of $W$, such that $w$ is $C^1$ on $W\bs S_w$, and 
$\H^{n-1}(S_w)=0$. 

$2^\circ$ A coflat $n-1$ form in $W$ is called closed if on $W\bs S_w$, $dw=0$; it is called a calibration, if it is closed, $||w||_\infty\le 1$ and $w(x)\ne 0, \forall x\in W\bs S_w$.

$3^\circ$ Let $w$ be a measurable $n-1$ form in $W$. Let $w_\perp$ be the vector field on $W$ defined as follows: if $w(x)=\sum_{i=1}^n f_i(x)e_1^*\wg\cdots e^*_{i-1}\wg e^*_{i+1}\wg\cdots \wg e^*_n$, then $w_\perp(x)=\sum_{i=1}^n f_i(x)e_i$. Here $f_i$ are measurable functions defined in $W$. It is easy to see that for a $C^1$ form $w$, $dw=0$ if and only if $w_\perp$ is of divergence zero.

$4^\circ$ Conversely, let $v$ be a measurable vector field in $W$ so that $v(x)\ne 0, \forall x\in W$. Let $v^\perp$ be the $C^1$ $n-1$ form in $W$ so that $(v^\perp)_\perp=v$. It is easy to see that for a $C^1$ vector field $v$, $div\ v=0$ if and only if $d(v^\perp)=0$.

$5^\circ$ A measurable vector field $v$ in $W$ is called a coflat vector field in $W$, if $v^\perp$ is a coflat $n-1$-form, and we denote by $S_v=S_{v^\perp}$ its singular set in $W$; it is called a divergence free coflat vector field, if $v^\perp$ is a closed coflat form.

\end{defn}

Note that every $n-1$-vector in $\R^n$ is simple, and the above defined $w_\perp$ is just a vector orthogonal to the $n-1$ space associated to $w$ at every point $x$.

\begin{defn}Let $d\le n$ be an integer. Let $E\subset \R^n$ be a closed set. We say that a $d$-integral current $T$ is associated to $E$, if $E=\spt (T)$, and $M(T)=\H^d(E)$ (in other words, the multiplicity is 1 almost everywhere).


\end{defn}

\begin{defn}
Let $U$ be an open subset of $\R^n$. Let $\beta$ be an $n-2$-integral current without boundary, such that $B=spt(\beta)\subset \pa U$. Let $E$ be a closed subset of $\bar U$ so that $E\cap \pa U=B$.

$1^\circ$ Let $w$ be a coflat calibration in $U$. We say that $E$ is calibrated by $w$ with respect to $\b$ in $U$, denoted by $E\in C(\b,w,U)$, if there exists an $n-1$ integral current $T$ associated to $E$, so that $\pa T=\b$, and $\H^{n-1}(E)=\lg T, w\rg$.

$2^\circ$ Let $v$ be a $C^1$ vector field in an open subset $W$ of $U$ so that $v(x)\ne 0,\forall x\in W$. Let $\cI$ be the class of all connected components of the integral curves of $v$ in $W$. We say that 

(1) $E$ is of multiplicity no more than 1 with respect to $v$, if for any $\gamma\in \cI$, $\gamma$ meets $E$ at at most one point;

(2) $E$ is of multiplicity no less than 1 with respect to $v$, if for every $x\in E\cap W$, there exists $\gamma\in \cI$ such that $x\in \gamma$;

(3) $E$ is of multiplicity 1 with respect to $v$, if it is of multiplicity no more than 1 and no less than one.

$3^\circ$ Let $w$ be a coflat calibration in $U$, let $W_w=U\bs S_w$. denote by $C_1(\b,w,U)=\{E\in C(\b,w,U): E$ is of multiplicity one with respect to $w_\perp|_{W_w}\}$. Also let $v(w)=w_\perp|_{W_w}$ be the $C^1$ vector field on $W_w$. 
\end{defn}

From Definition \tb{2.3} we know directly that 
\be C(\b,w,U)\subset \F_{ics}(n-1,B, \{\b\}, U).\ee

Let $U$ be an open set in $\R^n$, Let $v$ be a divergence free $C^1$ vector field in an open subset $W$ of $U$, so that $v(x)\ne 0, \forall x$. Let $E$ be a subset of $\bar U$ which is of multiplicity 1 with respect to $v$. Let $E_R=\{x\in E: \exists r>0$ so that $B(x,r)\cap E$ is a smooth manifold of dimension $n-1\}$ the regular part of $E$. Let $E_{R,v}=\{x\in E_R\cap W: v(x)\not\in T_xE_R\}\subset W\cap E_R$. Then $E_R$ and $E_{R,v}$ are both relatively open in $E$.

In the rest of the paper, whenever a coflat form $w$ in $U$ is involved, we always take $v=v(w)$ and $W=W_w$.

\begin{pro}Let $U$ be an open set in $\R^n$. Let $\b$ be an $n-2$-integral current without boundary, such that $B=spt(\b)\subset \pa U$. Let $w$ be a coflat calibration in $U$. If $E\in C(\b,w,U)$, then 

$1^\circ$ There exists a $n-1$-integral current $T$ associated to $E$ which is mass minimizing in $U\cup B$ and such that $\pa T=\b$. 

$2^\circ$ For $H^{n-1}$-almost all $x\in E$--in fact for all $x\in E_R\bs S_w$, $T_xE$ is the plane generated by $w(x)$, and $||w(x)||=1$. As a result, $E$ is of multiplicity no less than 1 with respect to $w_\perp$, and thus $C_1(\b,w,U)=\{E\in C(\b,w,U): E$ is of multiplicity no more than one with respect to $v(w)\}$. 

$3^\circ$ $E_R\bs S_w\subset E_{R,w_\perp}$, and $E\cap U\bs (S_w\cup E_{R,w_\perp})$ is of dimension at most $n-8$.

$4^\circ$ Let $[E_w]_{C^I_{n-1}}$ be the $n-1$-integral current $\H^{n-1}\lf_E\wg w_*$ (it is well an integral current by $2^\circ$). Then it is the unique $n-1$-integral current $T$ associated to $E$ so that $\pa T=\b$, and $\H^{n-1}(E)=\lg T, w\rg$. Moreover, it is mass minimizing in $U\cup B$. 
\end{pro}

\nd Since $E\in C(\b,w,U)$, by definition, there exists an $n-1$ integral current $T$ associated to $E$, so that $\pa T=\b$, and $\H^{n-1}(E)=\lg T, w\rg$. 

Take any $n$-current $S$ in $U\cup B$, then since $dw=0$ and $||w||\le 1$, we have
\be M(T+\pa S)\ge \lg T+\pa S, w\rg=\lg T,w\rg+\lg \pa S,w\rg=\H^{n-1}(E)+\lg S,dw\rg=M(T),\ee
which yields that $T$ is a mass minimizing current in $U\cup B$. This proves $1^\circ$.

Moreover, the above argument yields $M(T)=\lg T,w\rg$. Since $||w||\le 1$, we know that $T_xE$ is the plane generated by $w(x)$ and $||w(x)||=1$ for $\H^{n-1}$-almost all $x\in \spt(T)=E$. In particular, since $T_xE$ is continuous in $E_R\bs S_w$, and $E_R\bs S_w$ is relatively open, we know that for all $x\in E_R\bs S_w$, $T_xE$ is the plane generated by $w(x)$ and $||w(x)||=1$. This proves $2^\circ$. 

By $2^\circ$ and the definition we know that $E_R\bs S_w\subset E_{R,w_\perp}$.

 Since $T$ is a mass minimizing current of codimension 1 in $U\cup B$, by \cite{Fe70}, we know that $E\cap U=\spt(T)\bs \spt(\pa T)$ is a smooth submanifold except for a relatively closed set of Hausdorff dimension at most $n-8$. Thus $3^\circ$ is proved.
 
Now again let $T$ be any $n-1$-integral current $T$ associated to $E$ so that $\pa T=\b$, and $\H^{n-1}(E)=\lg T, w\rg$. Suppose that $T=\H^d\lf_E\wg w'$ where $w'$ is a measurable orientation on $E$ with $||w'||=1$. Then $\H^{n-1}(E)=\lg T, w\rg$ implies that $\lg w', w\rg=1$.  Since $||w'||=1$ and $||w||=1$,  and thus $w'=w_*$ for $\H^{n-1}$-a.e. $x\in E$. As a result, $T=[E_w]_{C^I_{n-1}}$. By the proof of $1^\circ$, it is also mass minimizing.
\qed

Next we will discuss deformation retracts along integral curves. (In smooth case, one can find discussions on retracts along curves that do not increase area in \cite{La91}.)

Let $U$ be an open set in $\R^n$. Let $v$ be a divergence free $C^1$ vector field in an open subset $W$ of $U$, so that $v(x)\ne 0, \forall x$. Let $E$ be a closed subset of $\bar U$ which is of multiplicity 1 with respect to $v$. Let $\theta:\cD\to W$ be the maximal flow of $v$, where $\cD$ is an open subset of $W\times \R$. For each $x\in W$, let $I_x=\{t\in \R: (x,t)\in \cD\}$.

Let $\gamma_x=\theta(\{x\}\times I_x)$. Set $\cD_0=\{(x,t)\in \cD: x\in E, t\in I_x\}$, and $W_0=W_0^E:=\theta(\cD_0)$. Then $E\cap W\subset W_0$. Define the map $p=p_v: W_0\to E\cap W$: for each $x\in W_0$, $p(x)\in E$ is such that $x\in \gamma_{p(x)}$. Since $E$ is of multiplicity 1 with respect to $v$, the map $p$ is well defined. And $p$ is surjective by definition. Note that $p(x)=x$ for $x\in E\cap W_0=E\cap W$.

\begin{defn}Let $U$ be an open set in $\R^n$. Let $v$ be a divergence free $C^1$ vector field in an open subset $W$ of $U$, so that $v(x)\ne 0, \forall x\in W$. Let $E$ be a closed subset of $\bar U$ which is of multiplicity 1 with respect to $v$. The above defined map $p_v: W_0\to E\cap W$ is called the projection along $v$ to $E$.

A closed set $F\subset \bar U$ is said to be of essentially full projection with respect to $(E,v)$ if $\H^{n-1}(E\bs p(F\cap W_0))=0$. The class of all sets of essentially full projection with respect to $(E,v)$ is denoted by $\C_p(E,v)$.
\end{defn}

It is easy to see that for any $x\in p^{-1}(E_{R,v})$, $p$ is $C^1$ in a neighborhood of $x$. As a result, $p$ is $C^1$ in the open set $p^{-1}(E_{R,v})$. 

The following lemma shows an "area decreasing property" for the map $p$. 

\begin{lem}Let $U$ be an open set in $\R^n$. Let $v$ be a divergence free $C^1$ vector field in an open subset $W$ of $U$, so that $v(x)\ne 0, \forall x\in W$. Let $E$ be a subset of $\bar U$ which is of multiplicity 1 with respect to $v$. Let $\theta,p=p_v, E_{R,v},W_0$ be as defined above. Then

$1^\circ$ Let $N\subset p^{-1}(E_{R,v})$ be an oriented compact smooth $n-1$-manifold with boundary, so that $p|_N$ is injective and $dp$ is non degenerate everywhere. Then we have 
\be\int_N v^\perp=\int_{p(N)} v^\perp\ee
where the orientation of $p(N)$ is induced by $p$.

$2^\circ$ Let $x\in p^{-1}(E_{R,v})$. Then for any $n-1$-vector $Q$ in $T_x (p^{-1}(E_{R,v}))=T_x\R^n\cong \R^n$, we have
\be \lg Dp(x)(Q), v^\perp\circ p(x)\rg=\lg Q, v^\perp(x)\rg.\ee

$3^\circ$ Let $k\le d$ be non negative integers. Let $\pi_1:\R^n\times \R^d\to \R^n$ and $\pi_2:\R^n\times\R^d\to \R^d$ be orthogonal projections. For each $x\in W_0\times \R^d$, let $x_i=\pi_i(x)$. Let $\varphi: W_0\times \R^d\to E\times \R^d: \varphi(x)=(p(x_1), x_2)$. 

Then for any $x\in p^{-1}(E_{R,v})\times \R^d$, and any $n-1+k$-vector $Q\in T_x(p^{-1}(E_{R,v})\times \R^d)=T_x(\R^n\times \R^d)\cong\R^n\times \R^d$, and any $k$-form $\xi$ on $\R^d$, we have
\be\lg D\varphi(x)(Q),((v^\perp\wg\xi)(\varphi(x))\rg=\lg Q, (v^\perp\wg \xi)(x)\rg,\ee
where $(v^\perp\wg\xi)(x):=v^\perp(x_1)\wg \xi(x_2)$.

$4^\circ$ Take all the notations as in $3^\circ$. Suppose in addition that $\H^{n-2}(E\cap W\bs E_R)=0$. 

Let $N\subset W_0\times \R^d$ be a $n-1+k$-rectifiable subset so that $\varphi|_N$ is injective, and let $u$ be a measurable orientation on $N$, that is, for any $x\in N$ so that the approximate tangent plane $T_xN$ exists, $u(x)\in \wg_{n-1}(T_xN)$ is a unit $n-1$-vector. Let $\varphi'$ be the restriction of $\varphi$ on $N$. Let $u'$ be the measurable orientation on $\varphi(N)$ associated to $u$, that is, $u'(y)=\frac{1}{||D\varphi'({\varphi'}^{-1}(y))||}D\varphi(u({\varphi'}^{-1}(y)))$ at any point $y$ so that $T_{{\varphi'}^{-1}(y)}N$ exists and $||D\varphi'({\varphi'}^{-1}(y))||\ne 0$.

Then for any $k$-form $\xi$ on $\R^d$
\be \lg\H^{n-1+k}\lf_{\varphi(N)}\wg u'(x), v^\perp\wg \xi\rg=\lg \H^{n-1+k}\lf_N\wg u, v^\perp\wg \xi\rg.\ee

In particular, in case $d=0$, we have
\be \lg\H^{n-1}\lf_{p(N)}\wg u'(x), v^\perp\rg=\lg \H^{n-1}\lf_N\wg u, v^\perp\rg,\ee
provided $N\subset W_0$ is $n-1$-rectifiable, $p|_N$ is injective, $u$ is a measurable orientation on $N$, and $u'$ is the orientation on $p(N)$ associated to $u$.

$5^\circ$ Let $w$ be a coflat $n-1$ form in $U$, and let $v=v(w)$ as defined as in Definition \tb{2.3}. Let $F\subset E$ be such that for some $M>0$, and a measurable orientation $u'$ on $F$ we have $\lg u', v\rg\ge M$ for $\H^{n-1}$-a.e. $x\in F$. Then for any $n-1$-rectifiable set $N\subset W_0$ with $p(N)\subset F$, we have
\be \H^{n-1}(p(N))\le \frac{||v|_{N}||_\infty}{M}\H^{n-1}(N).\ee
In particular, if $\b$ is an $n-2$-integral current without boundary, such that $B=spt(\b)\subset \pa U$, and $E\in C_1(\b,w,U)$, then for any $n-1$-rectifiable subset $N\subset W_0$, 
\be \H^{n-1}(p(N))\le \H^{n-1}(N).\ee
\end{lem}

\nd The conclusion and the proof of $1^\circ$-$3^\circ$ follows directly from \cite{product1} Lemma 3.6 $1^\circ$ -$3^\circ$.

The proof of $4^\circ$ and $5^\circ$ are slightly different from \cite{product1} Lemma 3.6 $4^\circ$ and $5^\circ$, but here we are facing an easier case: we do not have to consider the boundary of $E$:  

$4^\circ$ Take any $n-1+k$-rectifiable subset $N\subset W_0\times \R^d$ so that $\varphi|_N$ is injective. Let $u$ be a measurable orientation on $N$, that is, for any $x\in N$ so that the approximate tangent plane $T_xN$ exists, $u(x)\in \wg_{n-1}(T_xN)$ is a unit $n-1$-vector.
 
We first decompose $N$ as the disjoint union $N=N_1\cup N_2\cup N_3$, where $N_1=N\cap \varphi^{-1}[(E\cap W\bs E_R)\times\R^d]$, $N_2=\{x\in N\cap \varphi^{-1}((E_R\cap W)\times \R^d):$ the tangent plane $T_xN$ to $N$ at $x$ exists, and $||D\varphi(x)||=0\}$, and $N_3=N\cap \varphi^{-1}((E_R\cap W)\times \R^d)\bs N_2$. Then since $\varphi$ is injective, we know that $\varphi (N)$ is the disjoint union $\varphi(N_1)\cup \varphi(N_2)\cup \varphi(N_3)$. Hence it is enough to prove 
\be \int_{\varphi(N_i)}\lg u',v^\perp\wg\xi\rg d\H^{n-1+k}(x)=\int_{N_i}\lg u,v^\perp\wg\xi\rg d\H^{n-1+k}(x),i=1,2,3.\ee

For $i=1$: since $\H^{n-2}(E\cap W\bs E_R)=0$, we know that for $\H^{n-1+k}-a.e.x$ in $N_1$ (resp. $p(N_1)$), $\pi_1(T_xN_1)$ (resp. $\pi_1(T_x\varphi(N_1)$) is of dimension no more than $n-2$. This implies that $\lg u(x), (v^\perp\wg\xi)(x)\rg=0
$ for $\H^{n-1+k}-a.e.x$ in $N_1$, and $\lg u'(x), (v^\perp\wg\xi)(x)\rg=0
$ for $\H^{n-1+k}-a.e.x$ in $\varphi(N_1)$.
Therefore
\be\int_{\varphi(N_1)}\lg u',v^\perp\wg\xi\rg d\H^{n-1+k}(x)=\int_{N_1}\lg u,v^\perp\wg\xi\rg d\H^{n-1+k}(x)=0;\ee

For $i=2$, we know that $v^\perp$ is perpendicular to $T_xN_2$ at every point $x\in N_2$, and thus $\lg u(x), (v^\perp\wg\xi)(x)\rg=0$. On the other hand, by Sard's Theorem, $\H^{n-1+k}(\varphi(N_2))=0$. Hence we have
\be \int_{\varphi(N_2)}\lg u',v^\perp\wg\xi\rg d\H^{n-1+k}(x)=0=\int_{N_2}\lg u,v^\perp\wg\xi\rg d\H^{n-1+k}(x);\ee 

Finally for $i=3$, we apply \cite{Fe} Theorem 3.2.22, to the map $\varphi':N_3\to p(N_3)$ (where dim$N_3=$dim$p(N_3)$), and the function $\frac{1}{||D\varphi'(x)||}\lg u(x), (v^\perp\wg\xi)(x)\rg$ defined on $N_3$, and get
 \be \int_{\varphi(N_3)}\frac{1}{||D\varphi'({\varphi'}^{-1}(y))||}\lg u({\varphi'}^{-1}(y)), (v^\perp\wg\xi)({\varphi'}^{-1}(y))\rg d\H^{n-1+k}(y)=\int_{N_3} \lg u(x), (v^\perp\wg\xi)(x)\rg d\H^{n-1+k}(x).\ee
 By $3^\circ$, we know that $\lg u({\varphi'}^{-1}(y)), (v^\perp\wg\xi)({\varphi'}^{-1}(y))\rg=\lg D\varphi(u({\varphi'}^{-1}(y))), (v^\perp\wg\xi)(y)\rg$,  hence
\be\begin{split}&\frac{1}{||D\varphi'({\varphi'}^{-1}(y))||}\lg u({\varphi'}^{-1}(y)), (v^\perp\wg\xi)({\varphi'}^{-1}(y))\rg \\
&=\frac{1}{||D\varphi'({\varphi'}^{-1}(y))||}\lg D\varphi(u({\varphi'}^{-1}(y))), (v^\perp\wg\xi)(y)\rg\\
&=\lg \frac{1}{||D\varphi'({\varphi'}^{-1}(y))||}D\varphi(u({\varphi'}^{-1}(y))), (v^\perp\wg\xi)(y)\rg=\lg u'(y),(v^\perp\wg\xi)(y)\rg.
 \end{split}\ee
 Therefore by \tb{(2.13)} we get
 \be \int_{\varphi(N_3)}\lg u',v^\perp\wg\xi\rg d\H^{n-1+k}(x)=0=\int_{N_3}\lg u,v^\perp\wg\xi\rg d\H^{n-1+k}(x).\ee
 Suming up \tb{(2.10)-(2.12) and (2.15)}, we get \tb{(2.6)}.
 
 $5^\circ$ We decompose $N$ as the disjoint union of $N_i,1\le i\le 3$ as in $4^\circ$. Then we know that the $n-1$-Hausdorff measure of $p(N_i),i=1,2$ is zero.
 
 Now we take a measurable subset $N_4\subset N_3$, so that $p(N_4)=p(N_3)$, and the restriction of $p$ on $N_4$ is injective. 
 
By definition, $E\cap W\bs E_R\subset E\cap U\bs(S_w\cup E_{R,w_\perp})$, hence $\H^{n-2}(E\cap W\bs E_R)=0$ after Proposition \tb{2.4 $3^\circ$}. We can thus apply $4^\circ$. Again by Proposition \tb{2.4}, we know that for all $x\in p(N_4)$, $T_xp(N_4)$ is the $n-1$ subspace generated by $v^\perp$. Let $u$ be the orientation on $N_4$, so that $u'=(v^\perp)_*$ on $p(N_4)$. Then by \tb{$4^\circ$} we have
\be \begin{split}\H^{n-1}(p(N))&=\sum_{i=1}^3\H^{n-1}(p(N_i))=\H^{n-1}(p(N_3))=\H^{n-1}(p(N_4))\\
&\le \frac 1M\lg\H^{n-1}\lf{p(N_4)}\wg u',v^\perp\rg
=\frac 1M\lg\H^{n-1}\lf{N_4}\wg u,v^\perp\rg\\
&\le \frac{||v|_N||_\infty}{M}\H^{n-1}(N_4)\le \frac{||v||_\infty}{M}\H^{n-1}(N).\end{split}\ee
\qed

\begin{cor}Let $U$ be an open set in $\R^n$. Let $\b$ be an $n-2$-integral current without boundary, such that $B=spt(\b)\subset \pa U$. Let $w$ be a coflat calibration in $U$. If $E\in C_1(\b,w,U)$, then
$E$ minimizes the $n-1$-dimensional Hausdorff measure among all $n-1$-rectifiable sets with essentially full projection with respect to $(E,v(w))$, that is,
\be \H^{n-1}(E)=\inf_{F\in\C^{R,n-1}_p(E,v(w))}\H^{n-1}(F).\ee
\end{cor}

\nd We take all the notations as above. Let $p=p_{v(w)}$. 

Take any $F\in \C^{R,n-1}_p(E,v(w))$. 

by Lemma \tb{2.6 $5^\circ$}, we know that
\be \H^{n-1}(E)=\H^{n-1}(p(F\cap W_0))\le \H^{n-1}(F\cap W_0)\le \H^{n-1}(F).\ee\qed

Next let us study the minimalities for calibrated sets in $C_1(\b,w,U)$ among different classes. The key is to establish the essentially full projection property with respect to $(E,v(w))$. For this, we have to treat more carefully the singular set $S_w$ of $w$.

\begin{lem}Let $w$ be a coflat calibration in an open set $U\subset \R^n$. Let $W=W_w$, $v=v(w)$, and take all the notations as before. 

$1^\circ$ Set $E_S=\{x\in E\cap W: \bar \gamma_x\cap S_w\ne\emptyset\}$. Then $\H^{n-1}(E_S)=0$.

$2^\circ$  Similarly, set $E_B=\{x\in E\cap W: \bar \gamma_x\cap B\ne\emptyset\}$. Then $\H^{n-1}(E_B)=0$.\end{lem}

\nd For each $x\in E_S$, by definition there exists $y\in \bar \gamma_x\cap S_w$. Let $s(x)$ be ''the smallest distance in $\gamma_x$ between $x$ and $S_w$'', that is, 
\be s(x)=\inf\{\H^1(\gamma):\gamma\mbox{ is a connected subset of }\gamma_x, x\in \gamma\mbox{, and }\bar\gamma\cap S_w\ne\emptyset\}.\ee

Since $S_w$ is closed, $s(x)$ is never 0.

For each $k\ge 1$, let $E_S^k=\{x\in E_S, s(x)> \frac 1k\}\bs $. Then each $E_S^k\cap S_w=\emptyset$, and $E_S=\cup_{k\ge 1}E_S^k$. Thus it is enough to prove that $\H^{n-1}(E_S^k)=0$ for each $k\ge 1$.

Now fix any $k\ge 1$. 

Since $S_w$ is closed, for each $x\in E_S^k$, $d(x, S_w)>0$. Let $F_j=\{x\in E_S^k: d(x,S_w)>\frac 1j\}$. Then $E_S^k=\cup_{j\ge 1}F_j$. We would like to show that $\H^{n-1}(F_j)=0,\forall j$.

Fix $j\ge 1$. Take any $\e>0$. Since $\H^{n-1}(S_w)=0$, there exists countably many open balls $B_i,i\in\N$ of diameter less than $\min\{\frac 1k, \frac 1j\}$ so that $S_w\subset \cup_{i\in\N} B_i$, and $\sum_{i\in\N}(diam(B_i))^{n-1}<\frac{\e}{a_{n-1}}$, where $a_{n-1}$ is the $n-1$-area of the unit sphere in $\R^n$. As a result, if we let $\Sigma=\cup_{i\in\N}\pa B_i$, then $\Sigma$ is $n-1$-rectifiable and $\H^{n-1}(\Sigma)<\e$. 

By definition, for each $x\in F_j$, there exists a connected subset $\gamma\subset \gamma_x$, $\H^1(\gamma)>\frac 1k$, $x\in \gamma$ and $\bar\gamma\cap S_w\ne\emptyset$. Let $y\in \bar\gamma\cap S_w$. Take $i$ such that $y\in B_i$. Then since $d(x,S_w)>\frac 1j$ and $y\in S_w$, we know that $d(x,y)>\frac 1j$. Since $diam(B_i)<\frac 1j$, we know that $x\not\in B_i$. Now $\gamma$ is connected, which contains a point $x$ outside $B_i$ and a point inside $B_i$ (because $\bar\gamma\cap B_i\ne\emptyset$, hence $\gamma\cap B_i\ne\emptyset$). Take $z_x\in\gamma\cap \pa B_i$. Then by definition, $p(z_x)=x$.

The above argument asserts that for each $x\in F_j$, we can find $z_x\in \Sigma=\cup_i \pa B_i$ so that $p(z_x)=x$. Let $Z$ be the set of all $z_x, x\in F_j$, then $Z\subset W_0$, and $p(Z)=F_j$. Moreover since $Z\subset \Sigma$, we know that $Z$ is $n-1$-rectifiable and $\H^{n-1}(Z)<\e$.

Now by Lemma \tb{2.6 $5^\circ$}, we know that 
\be \H^{n-1}(F_j)=\H^{n-1}(p(Z))\le \H^{n-1}(Z)<\e.\ee

This holds for arbitrary $\e>0$, hence $\H^{n-1}(F_j)=0, \forall j\ge 1$, and therefore $\H^{n-1}(E_S^k)=0, \forall k\ge 1$. We have thus proved that $\H^{n-1}(E_S)=0$.

The proof of $2^\circ$ is exactly the same as that of $1^\circ$.\qed

Now let us discuss the important classes which belong to $\C_p(E,v(w))$ automatically. By Corollary \tb{2.7}, the set $E$ will minimize the Hausdorff measure among these classes.

\begin{pro}Let $U$ be an open bounded set in $\R^n$. Let $\b$ be an $n-2$-integral current without boundary, such that $B=spt(\b)\subset \pa U$.  Let $w$ be a coflat calibration in $U$. Let $E\in C_1(\b,w,U)$, then

$1^\circ$ $\F_{ihc}(E, n-1, B, \{[E_w]_{C^I_{n-1}}\}, U)\subset\C_p(E,v(w))$; 

$2^\circ$ $\F_{Al}(E,U)\subset \C_p(E,v(w))$;

$3^\circ$ If the $n-1$-dimensional integral rectifiable homology group of $U\cup B$ is trivial, then $\F_{ics}(n-1,B,\{\b\}, U)\subset \C_p(E,v(w))$.
\end{pro}

\nd $1^\circ$ Let $A\in \F_{ihc}(E, n-1, B,\{[E_w]_{C^I_{n-1}}\}, U)$. Set $T=[E_w]_{C^I_{n-1}}$ for short. Then there exists an integral $n-1$ current $S$ supported in $A$, and an integral $n$-current $R$ supported in $U\cup B$ so that $\pa R=T-S$. 

By Lemma \tb{2.8}, it is enough to prove that $E\cap W\bs (E_S\cup E_B)\subset p(spt(S)\cap W_0)$.

Fix any $x\in E\cap W\bs (E_S\cup E_B)$. Since $x\not\in E_S\cup E_B$, we know that $d(\bar\gamma_x, S_w\cup B)>0$, and $\bar\gamma_x\bs \gamma_x\subset\pa U$. Also since $spt(R)$ is compact, and $spt(R)\subset U\cup B$, we know that $d(\bar\gamma_x\cap \pa U, spt(R))>0$. Let $\d=\min\{d(\bar\gamma_x, S_w\cup B),d(\bar\gamma_x\cap \pa U, spt(R))\}>0$.

Let $\gamma_x^0=\gamma_x\bs B(\bar\gamma_x\cap \pa U, \frac \d2)$. Then $\gamma_x^0$ is compact, and hence $\inf_{z\in \gamma_x^0}w(z)>0$. Therefore $I_x^0:=\{t: \theta(x,t)\in \gamma_x^0\}$ is bounded, and hence compact. Let $a=\sup\{t\in I_x^0: [0,t]\subset I_x^0\}$, and $b=\inf\{t\in I_x^0: [t,0]\subset I_x^0\}$. Then $\theta(x,a),\theta(x,b)\in \pa B(\bar\gamma_x\cap \pa U, \frac \d2)$. 

By the continuity of $\theta$,  for any $t\in [a,b]$, there exists a neighborhood $V_t$ containing $\theta(x,t)$, and $\e_t>0$, so that 
$d(\theta(V_t\times (t-\e_t, t+\e_t)), \pa U\cup S_w)>\frac \d4$. Since $[a,b]$ is compact, there exists $\e>0$, such that $d(\theta(B(x,\e)\times [a,b]), \pa U\cup S_w)>\frac \d4$. Since $\gamma_x^0=\gamma_x\bs B(\bar\gamma_x\cap \pa U, \frac \d2)$, we can also ask that $\e$ is small enough so that, $E(x,\e):=B(x,\e)\cap E\subset E_R$, and $\theta(E(x,\e)\times \{a,b\})\subset B(\bar\gamma_x\cap \pa U, \frac{2\d}{3})$.

Now let us look at the "tube" region $A:=\theta(E(x,\e)\times (a,b))$. Let $S'=S\lf_A$, $T'=T\lf_A$, and $R'=R\lf_A$. Then we know that $\pa R'-(T'-S')$ is supported in $\pa A=\theta(\pa E(x,\e)\times [a,b])\cup \theta(E(x,\e)\times \{a,b\})$, where $\pa E(x,\e)$ is the manifold boundary of $E(x,\e)$. But we know that $\theta(E(x,\e)\times \{a,b\})\subset B(\bar\gamma_x\cap \pa U, \frac{2\d}{3})$, and $d(\bar\gamma_x\cap \pa U, spt(R))\}\ge \d$, therefore $\theta(E(x,\e)\times \{a,b\})\cap spt(R)=\emptyset$. As a result, $\pa R'-(T'-S')$ is supported in $\theta(\pa E(x,\e)\times [a,b])$, 
and thus $p_\sharp(\pa R')-p_\sharp(T'-S')$ is supported outside $E(x,\e)$. 

But $p_\sharp(R')$ is an $n$-integral current supported in the $n-1$ dimensional set $E$, hence $p_\sharp(R')=0$, and thus $\pa p_\sharp(R')=p_\sharp(\pa R')=0$. Since $p_\sharp(\pa R')-p_\sharp(T'-S')$ is supported outside $E(x,\e)$, we know that $p_\sharp(T'-S')$ is supported outside $E(x,\e)$. That is, $p_\sharp(S')\lf_{E(x,\e)}=p_\sharp(T')\lf_{E(x,\e)}=p_\sharp(T)\lf_{E(x,\e)}$, which is $\H^{n-1}\lf_{E(x,\e)}\wg w_*$ by definition. Hence $p_\sharp(S)\lf_{E(x,\e)}=p_\sharp(S')\lf_{E(x,\e)}=\H^{n-1}\lf_{E(x,\e)}\wg w_*$. This implies in particular that $E(x,\e)\subset p(W_0\cap \spt(S))$. Thus we get $1^\circ$.

$2^\circ$ By Proposition \tb{1.16}, we know that $F\in \F_{Al}(E,U)\subset\F_{ihc}(E, n-1, B, \{[E_w]_{C^I_{n-1}}\}, U)$. Hence $2^\circ$ follows from $1^\circ$ directly.

$3^\circ$ When the $n-1$-dimensional integral rectifiable homology group of $U\cup B$ is trivial, we know that $\F_{ics}(n-1,B,\{\b\}, U)=\F_{ihc}(E, n-1, B,\{[E_w]_{C^I_{n-1}}\}, U)$. Hence the conclusion follows directly from $2^\circ$.\qed


Combine Corollary \tb{2.7} Proposition \tb{2.9} and Proposition \tb{1.19}, we get

\begin{thm}Let $U$ be an open bounded set in $\R^n$. Let $\b$ be an $n-2$-integral current without boundary, such that $B=spt(\b)\subset \pa U$. Let $w$ be a coflat calibration in $U$. If $E\in C_1(\b,w,U)$, then

$1^\circ$ $E$ is a ''homological size minimizer'', that is: \be \H^{n-1}(E)=\inf_{F\in\F_{ihc}(E, n-1, B, \{[E_w]_{C^I_{n-1}}\}, U)}\H^{n-1}(F);\ee

$2^\circ$ $E$ is Almgren minimal of dimension $n-1$ in $U$;

$3^\circ$ If the $n-1$-dimensional integral rectifiable homology group of $U\cup B$ is trivial, then $E$ minimizes the $n-1$-Hausdorff measure among all $(n-1,B,\{\b\}, U)$-integral current spanners, that is,
\be \H^{n-1}(E)=\inf_{F\in\F_{ics}(n-1,B,\{\b\}, U)}\H^{n-1}(F).\ee
\end{thm}

\begin{rem}The multiplicity one property for $E$ is necessary. An example is given in the picture below: we have two horizontal circles with the same size whose centers are of the same first two coordinates in $\R^3$. We endow them with the same orientation, and let $B$ be their union. Let $E$ be the union of the two discs bounded by these two circles. Let $w=e_1^*\wg e_2^*$. Then it is a calibration, and $E\in C(\b,w,\R^3)$. But $E$ is obviously not of multiplicity one with respect to $w_\perp$, and hence when the distance between the two discs is relatively small, a pinching as in the picture will give a better competitor.

\centerline{\includegraphics[width=0.3\textwidth]{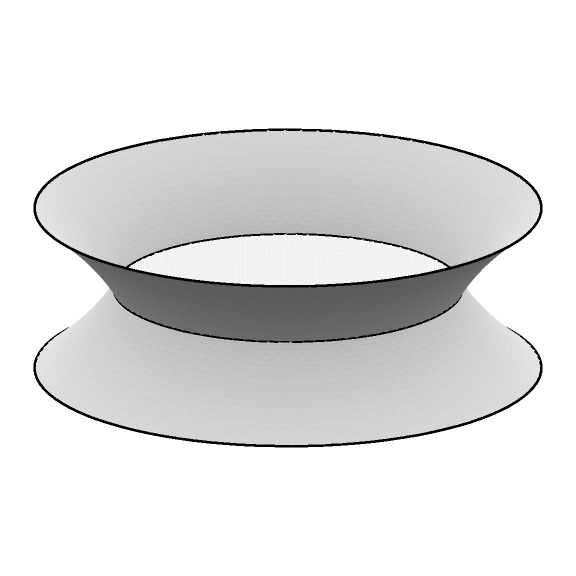} }
\end{rem}

\section{Minimality for paired calibrated sets}

Let $d\ge 2$ be an integer. Let $V\subset \R^d$ be a bounded strongly Lipschitz domain. Then the exterior unit normal vector $n(x)$ of $V$ exists for $\H^{d-1}$-a.e. $x\in \pa V$. 

For any $k\in \N$, set 
\be\begin{split} \cO(V,k)=\{\mo=(\O_1,\cdots \O_k): \O_i,1\le i\le k\mbox{ are mutually disjoint relatively open subsets of }\\
\pa V\mbox{ and each }\bar\O_i\mbox{ is a connected }d-1-\mbox{Lipschitz manifold with boundary}\}.
\end{split}\ee

For any $\mo=(\O_1,\cdots \O_k)\in \cO(V,k)$, let $B(\mo)=\pa V\bs (\cup_{i=1}^k \O_i)$. Then $\cup_{i=1}^k \pa \O_i\subset B(\mo)$ (here $\pa \O_i$ is the manifold boundary of $\bar\O_i$), but it may happen that $B(\mo)\bs (\cup_{i=1}^k \pa \O_i)\ne\emptyset$. 

We say that $\mo$ is a \textbf{$k$-Lipschitz partition} of $\pa V$, if $B(\mo)=\cup_{i=1}^k \pa \O_i$. In particular, $\H^{d-1}(B(\mo))=0$. Let $\cO_P(V,k)$ denote the set of all $k$-Lipschitz partitions of $\pa V$.
It is a subset of $\cO(V,k)$.

 Let $F\subset \bar V$ be a closed set, let $\mo=\{\O_1,\cdots \O_k\}$ be an element in $\cO(V,k)$. We say that $F$ \textbf{separates} $\mo$, if $F\cap \pa V=B(\mo)$, and $\O_i,1\le i\le k$ lie in different connected components of $\bar V\bs F$. We also say that $F$ is an \textbf{$\mo$-separator}. Denote by $\F_S(\mo)$ the class of all $\mo$-separators. 
 
 \medskip

For any fixed abelien group $G$, let $[\O_i]_G,1\le i\le k$ denote the element in $H_0(\pa V\bs B(\mo),G)$ represented by the connected component $\O_i$ with orientation pointing outward to $V$. Let $\pa\Omega_i$ also denote the manifold boundary of the oriented manifold $\Omega_i$. 
Set $[\mo]_G=\{[\O_i]_G-[\O_j]_G, 1\le i<j\le k\}\subset H_0(\pa V\bs B(\mo),G)$. And let $[\pa\mo]_{ic}$ denote the set of integral currents $\{[\pa \Omega_i]_{C^I_{d-2}}, 1\le i\le k\}$. 

We have the following relations between classes of sets:

\begin{lem}Let $V\subset \R^d$ ($d\ge 2$) be a bounded strongly Lipschitz domain. Let $\mo\in \cO(V,k)$ for some $k\in \N$. Then

$1^\circ$ 
\be \F_S(\mo)=\F_{tc}(d-1, G, B(\mo), [\mo]_G, V)\mbox{, for any non trivial abelian group }G;\ee
%

$2^\circ$ Let $E\subset \bar V$ be such that $E\cap \pa V=B(\mo)$, and for each $i$, there exists a $d-1$-integral current $\Gamma_i$ supported in $E$ so that $\Gamma_i$ is homologic to $[\O_i]_{C^I_{d-1}}$ in $V\cup \bar \O_i$. Then $E\in \F_S(\mo)$, and 
\be \F_{ihc}(E,d-1, B(\mo),\{\Gamma_i,1\le i\le k\} ,V)\subset \F_S(\mo).\ee

Conversely, \be\F_S^{R,d-1}(\mo)\subset \F_{ihc}(E,d-1, B(\mo),\{\Gamma_i,1\le i\le k\} ,V)\subset \F_{ics}(d-1, B(\mo), [\pa\mo]_{ic}, V).\ee

$3^\circ$ If the $d-1$-th integral rectifiable homology group $H^I_{d-1}(\bar V)$ is trivial, then
\be \F_{ics}(d-1, B(\mo), [\pa\mo]_{ic}, V)\subset \F_S(\mo).\ee
 \end{lem}
 
 \nd $1^\circ$ comes directly from definition;

$2^\circ$ Take any $E$ as in $2^\circ$. For each $i$, let $\Gamma_i$ be a $d-1$-integral current supported in $E$ such that there exists a $d$-integral current $\t D_i$ supported in $V\cup \bar \O_i$ so that $\pa \t D_i=\Gamma_i-[\O_i]_{C^I_{d-1}}$. Then there exists an integer valued density function $\theta$ in $D_i$ so that $\t D_i=\theta\H^d\lf_{D_i}\wg \omega$ where $\omega$ is an orientation $d$-form on $\R^d$ and $D_i$ is the support of $\t D_i$. 

Then $D_i$ is a set of finite perimeter, and $\bar \O_i\subset\pa D_i\subset \bar \O_i\cup E$. In particular, $\O_j\cap D_i=\emptyset$. Hence $\O_j\subset \bar V\bs D_i$. Now we have two disjoint relatively open subsets of $\bar V$: $\bar V\bs (D_i\cup E)$ and $D_i\bs E$, which contain $\O_j$ and $\O_i$ respectively, and whose union is $\bar V\bs E$, hence $\O_i$ and $\O_j$ belong to different connected components of $\bar V\bs E$. This holds for any $i\ne j$, hence $E\in \F_S(\mo)$.
 
Conversely, take any $F\in \F_S^{R,d-1}(\mo)$. By definition there exists a connected component $D_i$ of $\bar V\bs E$ so that $D_i\cap \pa V=\O_i$. Then $D_i$ is relatively open in $\bar V$, $\bar\O_i\subset\pa D_i\subset E\cup\bar\O_i$, and $\pa D_i$ is $d-1$ rectifiable with locally finite $\H^{d-1}$ measure, and by \cite{Fe} 4.5.11, we know that $D_i$ is a set of finite perimeter. Let $\t D_i$ denote the $d$-integral current induced by $D_i$, then $T_i:=\pa\t D_i-[\O_i]_{C^I_{d-1}}$ is a $d-1$-integral current supported in $F$. Thus $T_i$ is homologic to $[\O_i]_{C^I_{d-1}}$, and thus to $\Gamma_i$ in $V\cup\O_i$. Moreover since $spt(\Gamma_i-T_i)\subset V\cup \pa\O_i$, by constancy theorem we know that $spt(\t D_i)\subset V\cup \pa\O_i\subset V\cup spt(\Gamma_i)$. Hence $T_i$ is homologic to $\Gamma_i$ in $V\cup spt(\Gamma_i)$, $\forall 1\le i\le k$. That is, $F\in \F_{ihc}(E,d-1, B(\mo),\{\Gamma_i,1\le i\le k\} ,V)\subset \F_{ics}(d-1, B(\mo), [\pa\mo]_{ic}, V)$.

$3^\circ$ follows directly from $2^\circ$.\qed

\begin{lem}[cf. \cite{product1} Lemma 4.2]Let $V\subset \R^d$ ($d\ge 2$) be a bounded strongly Lipschitz domain. Let $\mo\in \cO_P(V,k)$ for some $k\in \N$. Then for each $F\in \F_S^{R,d-1}(\mo)$, we can find subsets $F_{i}, 1\le i\le  k$ of $F$, so that, if we denote by $F'$ the union of $F_i,1\le i\le k$, then

$1^\circ$ $V\bs F'$ has exactly $k$ connected components $D_i,1\le i\le k$, each $D_i$ is a set of finite perimeter, and we have the $d-1$-essentially disjoint union
\be \pa_* D_i=\Omega_i\cup F_i,1\le i\le k,\ee
where $\pa_*$ denotes the essential boundary. We also have the $d$-essentially disjoint union
\be \bar V=\cup_{1\le i\le k}\bar D_i;\ee

$2^\circ$ If we set $F_{ij}=F_i\cap F_j, 1\le i<j\le k$, then 
\be F'\mbox{ is the }d-1-\mbox{essentially disjoint union of }F_{ij},1\le i<j\le k.\ee
 \end{lem}

\nd See \cite{product1} Lemma 4.2.\qed

Next we introduce coflat paired calibrations.

\begin{defn}
$1^\circ$ For any open subset $V\subset \R^d$, let 
\be \cV_c^k(\bar V):=\{\nu=(v_1,\cdots v_k): v_i\mbox{ is a coflat divergence free vector field in a neighborhood of }V,1\le i\le k\}.\ee
For any $\nu\in \cV_c^k(\bar V)$, set $M(\nu)=\sup_{1\le i<j\le k}||v_i-v_j||_{L^\infty(\bar V)}$. Also let $S_\nu$ denote the union of the singular sets $S_{v_i}$: $S_\nu=\cup_{1\le i\le k}S_{v_i}$.

$2^\circ$ Let $V$ be a bounded strongly Lipschitz domain in $\R^d$, let $\nu\in \cV_c^k(\bar V)$, $\mo\in \cO_P(V,k)$, and $E\subset \bar V$ be closed. We say that $E$ is paired calibrated by $\nu$ with respect to $\mo$, denoted by $E\in PC(\mo, \nu)$, if

$1^\circ$ $E$ separates $\mo$;

$2^\circ$ $\H^{d-1}(E)=\frac{1}{M(\nu)}\sum_{i=1}^k\int_{\O_i}\lg n(x), v_i(x)\rg$.
\end{defn}

\begin{thm}Let $d\ge 2$ be an integer. Let $V$ be a bounded strongly Lipschitz domain in $\R^d$, so that $(\R^d,\pa V)$ admits a $C^1$ triangulation. Let $\nu\in \cV_c^k(\bar V)$ be such that $v_i-v_j$ does not vanish on $ V\bs S_\nu$ for any $i<j$. Let $\mo\in \cO_P(V,k)$. Let $E\in PC(\mo, \nu)$, then

$1^\circ$ $E$ minimizes the $d-1$-dimensional Hausdorff measure among all $\mo$-separators, that is, 
\be \H^{d-1}(E)=\inf_{F\in \F_S(\mo)}\H^{d-1}(F);\ee

$2^\circ$ $V\bs E$ has exactly $k$ connected components $\cC_i,1\le i\le k$, and $\bar\cC_i,1\le i\le k$ contain the $\O_i,1\le i\le k$ respectively; moreover, we have the $d-1$-essentially disjoint union $\pa \cC_i=\O_i\cup E_i$, where $E_i\subset E$ is a closed subset. In addition, $E=\cup_{1\le i\le k}E_i$ modulo null sets.

$3^\circ$ For $i\ne j$, define $E_{ij}:=E_i\cap E_j$. Then we have the $d-1$-essentially disjoint union $E=\cup_{1\le i<j\le k}E_{ij}$. Moreover, for each pair  $i<j$, for $\H^{d-1}$-almost all $x\in E_{ij}$, we have $v_i(x)-v_j(x)\perp T_xE$ and $||v_i(x)-v_j(x)||=M(\nu)$. 

$4^\circ$ Let $\cI_{ij}$ be the class of all integral curves of $v_i-v_j$ in $V\bs S_\nu$. For every $x\in E_{ij}\cap V\bs S_\nu$, let $\gamma^x\in \cI_{ij}$ be the element that contains $x$. Then for $\H^{d-1}$-almost all $x\in E_{ij}$, $\gamma^x\cap E_{ij}=\{x\}$, that is, it meets $E_{ij}$ at exactly one point. 
\end{thm}

\nd $1^\circ$ By Proposition \tb{1.19}, it suffices to prove that 
\be \H^{d-1}(E)=\inf_{F\in \F^{IR, d-1}_S(\mo)}\H^{d-1}(F).\ee

Fix any $F\in \F_S^{IR, d-1}(\mo)$.

We apply Lemma \tb{3.2} to $F$, and get $D_i$, $F_i$, $F'$ and $F_{ij}$ as in the conclusion of Lemma \tb{3.2}. Let $n_i(x)$ be the interior unit normal vector to $\pa_* D_i=\O_i\cup F_i$ at the point $x$, which exists for almost all $x\in \pa_* D_i$. Let $n(x)$ be the outer unit normal vector to $\pa V$ at $x\in \pa V$, which exists for $\H^{d-1}$-almost all $x\in \pa V$. Then $n_i(x)=-n(x)$ for $\H^{d-1}$-almost all $x\in \O_i$.

Now for each $1\le i\le k$, we apply the divergence theorem to the vector field $v_i$ and the region $D_i$, and get
\be \int_{\pa_* D_i}\lg n_i(x), v_i(x)\rg d\H^{d-1}(x)=0.\ee
Since $n_i(x)=-n(x)$ for $\H^{d-1}$-a.e. $x\in \O_i$. Hence we have
\be\begin{split}0&=\int_{\pa_* D_i}\lg n_i(x), v_i(x)\rg d\H^{d-1}(x)=\int_{\O_i}\lg n_i(x), v_i(x)\rg d\H^{d-1}(x)+\int_{F_i}\lg n_i(x), v_i(x)\rg d\H^{d-1}(x)\\
&=\int_{F_i}\lg n_i(x), v_i(x)\rg d\H^{d-1}(x)-\int_{\O_i}\lg n(x), v_i(x)\rg d\H^{d-1}(x).\end{split}\ee
We sum over $i$, and get
\be \begin{split}\sum_{1\le i\le k}\int_{\O_i}\lg n(x), v_i(x)\rg d\H^{d-1}(x)=\sum_{1\le i\le k}\int_{F_i}\lg n_i(x), v_i(x)\rg d\H^{d-1}(x).\end{split}\ee

Note that for almost all $x\in F_{ij}$, $n_i(x)=-n_j(x)$, hence by Lemma \tb{3.2 $2^\circ$},
\be\begin{split}&\sum_{1\le i\le k}\int_{F_i}\lg n_i(x), v_i(x)\rg d\H^{d-1}(x)=\sum_{1\le i<j\le k}\int_{F_{ij}}\lg n_i(x), v_i(x)-v_j(x)\rg d\H^{d-1}(x)\\
&\le \sum_{1\le i<j\le k}\int_{F_{ij}}||v_i(x)-v_j(x)|| d\H^{d-1}(x)\le  \sum_{1\le i<j\le k} M(\nu)\H^{d-1}(F_{ij})\\
&=M(\nu)\H^{d-1}(F')\le M(\nu)\H^{d-1}(F).
\end{split}\ee
Combine with \tb{(3.14)} we get
\be \sum_{1\le i\le k}\int_{\O_i}\lg n(x), v_i(x)\rg d\H^{d-1}(x)\le M(\nu)\H^{d-1}(F')\le M(\nu)\H^{d-1}(F).\ee
Since $E\in PC(\mo,\nu)$, we know that
\be  M(\nu)\H^{d-1}(E)=\sum_{1\le i\le k}\int_{\O_i}\lg n(x), v_i(x)\rg d\H^{d-1}(x)\le M(\nu)\H^{d-1}(F),\ee
for all $F\in \F_S^{IR,d-1}(\mo)$. Thus we get \tb{(3.11)}. And hence $1^\circ$ is proved.

$2^\circ$ By $1^\circ$, Lemma \tb{3.1 $2^\circ$} and Proposition \tb{1.16}, we know that $E\cap V$ is an Almgren minimal set of dimension $d-1$ in $V$. As a result, $E$ is $d-1$-rectifiable, i.e., $E\in \F_S^{R,d-1}(\mo)$.

Since $E\in \F_S^{R,d-1}(\mo)$, we define $D_i=D_i(E)$ and the 
$E_i,1\le i\le k$, $E_{ij}, 1\le i<j\le k$ and $E'=\cup_{1\le i<j\le k}E_{ij}$ as in Lemma \tb{3.2}. Then all the inequalities in \tb{(3.15) and (3.16)} are in fact equalities for $F=E$. That is:
\be \lg n_i(x), v_i(x)-v_j(x)\rg=||v_i(x)-v_j(x)||=M(\nu), \H^{d-1}-a.e. x\in E_{ij},\ee
and 
\be \H^{d-1}(E\bs E')=0.\ee

In particular, \tb{(3.18) and (3.19)} tells that for $\H^{d-1}$-a.e. $x\in E$, there exists $1\le i<j\le k$ so that $x\in E_{ij}\subset \pa_* D_i\cap \pa_* D_j$.
This implies that $D_1,\cdots, D_k$ are all the connected components of $V\bs E$, for if there exists another one, then its boundary will be of positive $\H^{d-1}$ measure, and will not belong to any of these $E_{ij}$. Thus we get $2^\circ$.

$3^\circ$ follows directly from \tb{(3.18)}.

$4^\circ$ Fix any $i\ne j$. Since $v_i-v_j$ never vanishes on $V\bs S_\nu$, the integral curves of $v_i-v_j$ cover all $V\bs S_\nu$, so that each $x\in E\cap V\bs S_\nu$ is passed by a $\gamma^x$.

Let $\gamma\in \cI_{ij}$. Then $\gamma\subset V$. 

By \tb{$3^\circ$}, for $\H^{d-1}$-a.e. $x\in E_{ij}$, and any $y\in E_{ij}\cap\gamma^x$, $v_i(y)-v_j(y)\not\in T_yE_{ij}$, hence $\lg n_i(y), v_i(y)-v_j(y)\rg\ne 0$. Note that $\lg n_i(y), v_i(y)-v_j(y)\rg>0$ means that the directed curve $\gamma^x$ (with orientation $v_i(y)-v_j(y)$) goes from $\cC_j$ to $\cC_i$ through $y$, and $\lg n_i(y), v_i(y)-v_j(y)\rg<0$ the converse. Since $\gamma^x$ is connected, then when we go along $\gamma^x$ towards the direction of $v_i$, we will pass the boundary $\pa\cC_i\cap \pa\cC_j$ through the above two types of points in turn. In particular, if $\gamma^x\cap E_{ij}=\gamma^x\cap (\pa\cC_i\cap \pa\cC_j)$ contains more than 2 points, then there exists a point $y\in \gamma^x\cap E_{ij}$ through which $\gamma^x$ goes from $\cC_j$ to $\cC_i$, i.e. $\lg n_i(y), v_i(y)-v_j(y)\rg<0$. This contradicts the fact that $\lg n_i(y), v_i(y)-v_j(y)\rg=M(\nu)>0$.

Thus we get $4^\circ$.
\qed

After Theorem \tb{3.4}, let $[E_i]_{C^I_{d-1}}$ be the integral current associated to $E_i$, so that $\pa[\cC_i]_{C^I_d}=[E_i]_{C^I_{d-1}}-[\O_i]_{C^I_{d-1}}$.

\begin{thm}Let $d\ge 2$ be integer. Let $V$ be a bounded strongly Lipschitz domain in $\R^d$, so that $(\R^d,\pa V)$ admits a $C^1$ triangulation. Let $\nu\in \cV_c^k(\bar V)$ be such that $v_i-v_j$ does not vanish on $V\bs S_\nu$ for any $i<j$. Let $\mo\in \cO_P(V,k)$. Let $E\in PC(\mo, \nu)$, then the following assertions hold:

$1^\circ$ $E$ is a "homological minimizer", that is: 
\be \H^{d-1}(E)=\inf_{F\in\F_{ihc}(E,d-1, B(\mo),\{[E_i]_{C^I_{d-1}},1\le i\le k\} ,V)}\H^{d-1}(F).\ee

$2^\circ$ $E$ is a $d-1$-dimensional Almgren minimal set in $V$;

$3^\circ$ If the $d-1$-th integral rectifiable homology group $H^I_{d-1}(\bar V)$ is trivial, then
\be \H^{d-1}(E)=\inf_{F\in\F_{ics}(d-1, B(\mo),[\pa\mo]_{ic}, V)}\H^{d-1}(F).\ee

\end{thm}

\nd  It follows directly from Theorem \tb{3.4}, Lemma \tb{3.1 $2^\circ$, $3^\circ$} and Proposition \tb{1.16}.
\qed

\begin{rem}In the theorem, we ask a relatively high regularity for the domain $V$. But in many cases, even when $V$ is not regular enough, we can take a regular domain $V'$ which is slightly larger than $V$, and extend $v_i, 1\le i\le k$ to $V'$ to apply the theorem, so that the paired calibrated set is minimal in $V'$, and thus in $V$. Also, if $V$ contains the convex hull of $B(\mo)$, and every point of $B(\mo)$ belongs to the boundary of its convex hull, then we can still get the minimality of any set in $PC(\mo,\nu)$ due to the regularity of any convex domain and the convex hull property of any minimizer.
\end{rem}

\section{Minimality of the product}

Now we are going to discuss the various minimality results for the product of a codimensional 1 calibrated set of multiplcity 1 and a paired calibrated set.

Let $n\ge 2,d\ge 2$ be integers. Let $U$ be an open bounded set in $\R^n$, and $V$ be a bounded strongly Lipschitz domain in $\R^d$, so that $(\R^d,\pa V)$ admits a $C^1$ triangulation. Let $\b$ be an $n-2$-integral current without boundary, such that $B=spt(\b)\subset \pa U$, and let $\mo=(\O_1,\cdots,\O_k)\in \cO_P(V,k)$. Let $E\in \F_{ics}(n-1,B,\{\b\},U)$, and $F\in \F_{ics}(d-1, B(\mo), [\pa\mo]_{ic}, V)$. Then $U\times V$ is open in $\R^{n+d}$, and $\pa (U\times V)=(\pa U\times V)\cup (U\times \pa V)$. The set $E\times F\subset \overline{U\times V}$, and $(E\times F)\cap \pa(U\times V)=(B\times F)\cup (E\times B(\mo))$. 

Since $E\in \F_{ics}(n-1,B,\{\b\},U)$, there exists a $n-1$-integral current $T$ supported in $E$ so that $\pa T=\b$; similarly, since $F\in \F_{ics}(d-1, B(\mo), [\pa\mo]_{ic}, V)$, for each $1\le i\le k$, there exists a $d-1$-integral current $T_i$ supported in $F$ so that $\pa T_i=[\pa \O_i]_{C^I_{d-2}}$. As a result, the $n-d-2$-integral current $T\times T_i$ satisfies that 
\be\pa (T\times T_i)=(\pa T)\times T_i+(-1)^n T\times(\pa T_i)=\b\times T_i+(-1)^n T\times [\pa \O_i]_{C^I_{d-2}}.\ee

Now if $E\in C_1(\b,w,U)$ for a coflat calibration $w$ in $U$, then we take $T=[E_w]_{C^I_{n-1}}=\H^{n-1}\lf_E\wg w_*$ as in Proposition \tb{2.4}. 

In addition if $F\in PC(\mo,\nu)$ for some $\nu\in \cV_c^k(\bar V)$ such that $v_i-v_j$ does not vanish on $V$ for any $i<j$. Let $\cC_i, F_i,1\le i\le k$ be defined by \tb{$2^\circ$} in Theorem \tb{3.4}. Let $[F_i]_{C^I_{d-1}}$ be the integral current associated to $F_i$, so that $\pa[\cC_i]_{C^I_d}=[F_i]_{C^I_{d-1}}-[\O_i]_{C^I_{d-1}}$. We then take $T_i=[F_i]_{C^I_{d-1}}$ in \tb{(4.1)}. Thus \tb{(4.1)} becomes
\be \pa ([E_w]_{C^I_{n-1}}\times [F_i]_{C^I_{d-1}})=\b\times [F_i]_{C^I_{d-1}}+(-1)^{n-1} [E_w]_{C^I_{n-1}}\times [\pa \O_i]_{C^I_{d-2}}.\ee

For $1\le i\le k$, let 
\be\begin{split}[\cU_i]_{C^I_{n+d-2}}&=\pa ([E_w]_{C^I_{n-1}}\times [\cC_i]_{C^I_d})+(-1)^n[E_w]_{C^I_{n-1}}\times [F_i]_{C^I_{d-1}}\\
&=\b\times [\cC_i]_{C^I_d}+(-1)^{n-1}[E_w]_{C^I_{n-1}}\times\pa  [\cC_i]_{C^I_d}+(-1)^n[E_w]_{C^I_{n-1}}\times [F_i]_{C^I_{d-1}}\\
&=\b\times [\cC_i]_{C^I_d}+(-1)^{n-1}[E_w]_{C^I_{n-1}}\times\{[F_i]_{C^I_{d-1}}-[\O_i]_{C^I_{d-1}}\}+(-1)^n[E_w]_{C^I_{n-1}}\times [F_i]_{C^I_{d-1}}\\
&=\b\times [\cC_i]_{C^I_d}+(-1)^n[E_w]_{C^I_{n-1}}\times[\O_i]_{C^I_{d-1}}.\end{split}\ee

Then by definition, we know that $[\cU_i]_{C^I_{n+d-2}}$ is homologic to $(-1)^{n-1}[E_w]_{C^I_{n-1}}\times [F_i]_{C^I_{d-1}}$ in $\bar U\times \bar V$. Let $\cU_i=spt([\cU_i]_{C^I_{n+d-2}})$. Then $\cU_i\subset \pa (U\times V)$.

%

The first main theorem of the section is the following:

\begin{thm}Let $n\ge 2,d\ge 2$ be integers. 

Let $U$ be an open bounded set in $\R^n$. Let $\b$ be an $n-2$-integral current without boundary, such that $B=spt(\b)\subset \pa U$.  Let $E\in C_1(\b,w,U)$, where $w$ is a coflat calibration in $U$, such that its singular set $S_w\subset E$, and $\H^{n-2}(S_w)=0$. 

Let $V\subset \R^d$ be a bounded strongly Lipschitz domain, so that $(\R^d,\pa V)$ admits a $C^1$ triangulation. Let $F\in PC(\mo,\nu)$, where $\mo=(\Omega_1,\cdots, \Omega_k)\in \cO_P(V,k)$ for some $k\in \N$, and $\nu\in \cV_c^k(\bar V)$ is such that $v_i-v_j$ does not vanish on $V\bs S_\nu$ for any $i<j$.

Suppose that $U\times V$ is a SR-domain. 

Let $H=\{[E_w]_{C^I_{n-1}}\times [F_i]_{C^I_{d-1}}, 1\le i\le k\}$, where $F_i$ and $[F_i]_{C^I_{d-1}}$ are as defined above. 

Then  \be \H^{n+d-2}(E\times F)=\inf\{\H^{n+d-2}(A):A\in \F_{ihc}(E\times F, n+d-2,(B\times F)\cup (E\times B(\mo)),H,U\times V)\}.\ee
\end{thm}

The proof of Theorem \tb{4.1} will keep us busy for the rest of this section. Before proving the theorem, we will first introduce some notations, and then give two propositions which asserts that we can restrict our attention to competitors with better properties. Finally we finish our proof.

 Let $\Gamma=(B\times F)\cup (E\times B(\mo))=(E\times F)\cap \pa (U\times V)$. Let $\Gamma_i=(B\times F_i)\cup (E\times \pa\O_i)$.
 Note that $U\times V$ is a bounded strongly Lipschitz domain in $\R^{n+d}$, by Proposition \tb{1.19} it is enough to prove that 
 \be \H^{n+d-2}(E\times F)\le \inf\{\H^{n+d-2}(A): A\in \F^{IR,d}_{ihc}(E\times F, n+d-2,\Gamma,H,U\times V)\}.\ee
 
 Let us first introduce some notations. 
 
 Let $W=W_w$, $v=v(w)$ (cf. Definition \tb{2.3}), and take the notations $W_0, p=p_v, \theta, \gamma_x$ etc. as before. 
 
Let $E'=E\times \bar V$. Let $\zeta$ be a unit $d$-vector of $\R^d$, and let $w'=w\wg \zeta^*$, then it is a coflat calibration in $U\times V$, and its singular set is $S_{w'}=S_w\times V$. Let $W'=W_{w'}=W\times V\subset U\times V$, and let $v'=v(w')$ in $U\times V$ (cf. Definition \tb{2.3}). Also, let $W'_0\subset U\times V$, $\theta'$, $\gamma'_z$ for $z\in W'$, $p':W'_0\to E'\bs S_{w'}$ be defined similarly. We extend $p'$ to $D':=E'\cup W'_0$ by defining $p'(z)=z$ for all $z\in E'$.
 
 Let $B'=E'\cap \pa(U\times V)$. Let $\b'=\pa [E'_{w'}]_{C^I_{n+d-1}}$ where $[E'_{w'}]_{C^I_{n+d-1}}=\H^{n+d-1}\lf_{E'}\wg w'_*$.
  
Note that $v'(x,y)=(v(x),0)$ for all $(x,y)\in U\times V$; $\theta'(z,t)=(\theta(x,t),y)$, $\gamma'_z=\gamma_x\times \{y\}$ for $z=(x,y)\in W'$; $W'_0=W_0\times V$, $p'(x,y)=(p(x), y)$ for all $(x,y)\in W'_0$. The set $E'\in C_1(\b', w', U\times V)$. By Theorem \tb{2.10}, $E'$ is an $n+d-1$-Almgren minimal set in $U\times V$. 
 
 
 Set $\cT_\e=B(\pa(U\times V),\e)=\{z\in U\times V: d(z, \pa(U\times V))<\e\}$. 
  %

After introducing the above notations, let us give two simplifications: Proposition \tb{4.2} and Proposition \tb{4.3}. Proposition \tb{4.2} says that we can only consider competitors which coincides with $E'=E\times \bar V$ near $\pa (U\times V)$. Based on this, Proposition \tb{4.3} says that modulo adding a set of arbitrary small $\H^{n+d-2}$-measure, the projection of such a competitor under $p'$ can "separate" the regions $\cU_i, 1\le i\le k$ (which form a "partition" of the boundary of $E\times \bar V$) inside the set $E'=E\times \bar V$. Recall that $[\cU_i]_{n+d-2}^I$ is homologic to $[E_w]_{C^I_{n-1}}\times [F_i]_{C^I_{d-1}}$, $1\le i\le k$.

 \begin{pro}Under the hypotheses of Theorem \tb{4.1}: For each $\e>0$, and for any $A\in \F_{ihc}^{IR,n+d-2}(E\times F, n+d-2,\Gamma,H,U\times V)\}$, there exists $\d>0$ and $A'\in \F_{ihc}^{IR,n+d-2}(E\times F, n+d-2,\Gamma,H,U\times V)\}$, so that $A'\cap \cT_{4\d}\subset E'$,  $\H^{n+d-2}(A')<\H^{n+d-2}(A)+\e$.
 \end{pro}
 
 \nd Recall that $U\times V$ is a SR-domain. Let $L, \e_0$ be its SR-constants. Take any $\eta<\e_0$. Let $\varphi:\overline{U\times V}\times [0,1]$ be a $(L,\eta)$-self retract of $U\times V$. Let $\Theta_\eta=(E\times F)\bs (U\times V)^-_{2\eta}=(E\times F)\cap \cT_{2\eta}$. Let $A_1=A\cup \Theta_\eta$.  Set $A'=\varphi_1(A_1)\cup \Theta_\eta$. Since $\varphi_1(A_1)\subset U\times V$, there exists $\d=\d_\eta\in (0, \eta/4)$ so that $\varphi_1(A_1)\subset (U\times V)^-_{4\d}$. Thus we know that $A'\cap \cT_{4\d}=\Theta_\eta\subset E\times F\subset E'$. Moreover, 
 \be \begin{split}\H^{n+d-2}(A')&\le\H^{n+d-2}[\varphi_1(A_1)]+\H^{n+d-2}(\Theta_\eta)\\
 &\le \H^{n+d-2}[\varphi_1(A_1\cap (U\times V)^-_{2\eta})]+\H^{n+d-2}[\varphi_1(A_1\cap \cT_{2\eta}]+\H^{n+d-2}(\Theta_\eta)\\
 &\le \H^{n+d-2}(A\cap (U\times V)^-_{2\eta})+L^{n+d-2}\H^{n+d-2}(A_1\cap  \cT_{2\eta})+\H^{n+d-2}(\Theta_\eta)\\
 &\le \H^{n+d-2}(A)+L^{n+d-2}[\H^{n+d-2}(A\cap \cT_{2\eta})+\H^{n+d-2}(\Theta_\eta)]+\H^{n+d-2}(\Theta_\eta)\\
 &=\H^{n+d-2}(A)+L^{n+d-2}\H^{n+d-2}(A\cap \cT_{2\eta})+(L^{n+d-2}+1)\H^{n+d-2}(\Theta_\eta).
 \end{split}\ee
 
 Now since $\H^{n+d-2}[A\cap \pa(U\times V)]=\H^{n+d-2}[(E\times F)\cap \pa(U\times V)]=0$, we know that
 \be\lim_{\eta\to 0} L^{n+d-2}\H^{n+d-2}(A\cap \cT_{2\eta})+(L^{n+d-2}+1)\H^{n+d-2}(\Theta_\eta)=0.\ee
 Hence there exists $\eta$ so that $L^{n+d-2}\H^{n+d-2}(A\cap \cT_{2\eta})+(L^{n+d-2}+1)\H^{n+d-2}(\Theta_\eta)<\e$, that is
 \be \H^{n+d-2}(A')<\H^{n+d-2}(A)+\e.\ee
 
 Fix such an $\eta$, and thus a $\d=\d_\eta$. The rest is to prove that the corresponding $A'$ satisfies that $A'\in \F_{ihc}^{IR,n+d-2}(E\times F, n+d-2,\Gamma,H,U\times V)\}$ for each $\eta$. Since $A'$ is the union of Lipschitz images of $A$ and $\Theta_\eta$, both of which are $n+d-2$-integral regular, we know that so is $A'$. 
 
 Now let $1\le i\le k$. Let $S'_i=[E_w]_{C^I_{n-1}}\times [F_i]_{C^I_{d-1}}$ for short. Since $A\in \F_{ihc}(E\times F, n+d-2,\Gamma,H,U\times V)\}$, there exists an $n+d-2$ integral current $S_i$ supported in $A$, and a $n+d-1$-integral current $R_i$ supported in $(U\times V)\cup \Gamma_i$ so that $\pa R_i=S_i-
S'_i$. Then we have
\be \pa R_i=S_i-S'_i\lf_{\cT_{2\eta}}-S'_i\lf_{(U\times V)^-_{2\eta}}.\ee
 Now let $S_i^1=S_i-S'_i\lf_{\cT_{2\eta}}$, and $S_i^2=S'_i\lf_{(U\times V)^-_{2\eta}}$. Then they are homologic, $spt(S_i^1)\subset A_1$, and $spt(S_i^2)\subset (U\times V)^-_{2\eta}$.
 
 Next, since $\varphi$ is a Lipschitz deformation, we know that $\varphi_{1,\s}(S_i^1)$ is homologic to $\varphi_{1,\s}(S_i^2)$ in $(U\times V)^-_{4\d}$. Since $spt(S_i^2)\subset (U\times V)^-_{2\eta}$, and $\varphi_1|_{(U\times V)^-_{2\eta}}=id$, we know that $\varphi_{1,\s}(S_i^2)=S_i^2$, and hence $\varphi_{1,\s}(S_i^1)\sim_{H^I_{n+d-2}((U\times V)^-_{2\eta})}S_i^2$. Then since $S'_i\sim_{H^I_{n+d-2}(E\times F_i)} S'_i=S'_i\lf_{\cT_{2\eta}}+S_i^2$, which means that $S_i^2\sim_{H^I_{n+d-2}(E\times F_i)}S'_i-S'_i\lf_{\cT_{2\eta}}$, we have
 \be \varphi_{1,\s}(S_i^1)\sim_{H^I_{n+d-2}((U\times V)^-_{2\eta}\cup (E\times F_i))}S'_i-S'_i\lf_{\cT_{2\eta}},\ee
 and hence
 \be \varphi_{1,\s}(S_i^1)+S'_i\lf_{\cT_{2\eta}}\sim_{H^I_{n+d-2}((U\times V)\cup [(E\times F_i)\cap \pa (U\times V)])}S'_i.\ee
 Note that $spt(\varphi_{1,\s}(S_i^1)+S'_i\lf_{\cT_{2\eta}})\subset A'$. We have thus proved that $A'\in \F_{ihc}^{IR,n+d-2}(E\times F, n+d-2,\Gamma,H,U\times V)\}$.\qed
 
 Note that in the proof of Proposition \tb{4.2} we do not need the hypothesis $\H^{n-2}(S_w)=0$. But it will be needed in the following proposition.
  
  \begin{pro}Under the hypotheses of Theorem \tb{4.1}. Let $A\in \F_{ihc}^{IR,d}(E\times F, n+d-2,\Gamma,H,U\times V)\}$ so that there exists $\d>0$ such that $A\cap \cT_{4\d}\subset E'$. Then for any $\e>0$, and any $1\le i\le k$, there exists a $n+d-1$ integral current $Q_i$ supported on $E'$
so that $\cU_i\subset spt(\pa Q_i))$, $spt(\pa Q_i-[\cU_i]_{C^I_{n+d-2}})$ is $\H^{n+d-2}$ essentially disjoint from $\pa(U\times V)$, and $\H^{n+d-2}(spt(\pa Q_i)\bs [p'(A\cap D')\cup \cU_i])<\e$.
 \end{pro}
 
\nd Since $A\cap \cT_{4\d}\subset E'$, for each $1\le i\le k$, there exists a $n+d-1$ integral current $R_i$ supported in $(U\times V)\cup \Gamma_i$, so that $spt(\pa R_i-[E_w]_{C^I_{n-1}}\times [F_i]_{C^I_{d-1}})\subset A$, and $spt(R_i)\cap \cT_{4\d}\subset E'$.

Fix any $1\le i\le k$. Let $T_i=\pa R_i-[E_w]_{C^I_{n-1}}\times [F_i]_{C^I_{d-1}}$. Then $spt(T_i)\subset A$. Let $R_i^1=R_i+(-1)^n([E_w]_{C^I_{n-1}}\times [\cC_i]_{C^I_d})$, then $\pa R_i^1=T_i+(-1)^n[\cU_i]_{C^I_{n+d-2}}$. Hence $spt(\pa R_i^1-(-1)^n[\cU_i]_{C^I_{n+d-2}})\subset A$. Moreover since $spt([E_w]_{C^I_{n-1}}\times [\cC_i]_{C^I_d})\subset E'$, we still have $spt(R^1_i)\cap \cT_{4\d}\subset E'$.

We would like to "project" $R_i^1$ to $E'$ along $p'$. So we have to deal with the singular set $S'_w$.
 
Since $E'$ is $n+d-1$-Almgren minimal in $U\times V$, by the monotonicity of density and the Ahlfors regularity for reduced minimal sets (cf. \cite{DJT}), there exists a constant $C$ that depends only on $n$, so that for each $z\in E'\cap(U\times V)$, the density ratio function $r\mapsto\theta_z(r)=\frac {\H^{n+d-1}(E'\cap B(z,r)}{r^{n+d-1}})$ is non decreasing on $(0, \frac 12d(z, \pa (U\times V)))$ and bounded by $C$, thus the $\H^{n+d-1}$-density $\theta(z)=\lim_{r\to 0}\frac {\H^{n+d-1}(E'\cap B(z,r)}{r^{n+d-1}}$ of $E'$ on $z$  exists and is no more than $C$. Therefore for each $z\in S'_w\subset E'$, there exists $r_z>0$ so that $B(z,r_z)\subset U\times V$, and for any $r<r_z$, the density ratio $\theta_z(r)\in [\theta(z), \frac32\theta(z))$. This means, for any $0<t<s<r_z$, 
 \be\theta(z)s^{n+d-1}\le \H^{n+d-1}(B(z,s)\cap E')\le \frac32\theta(z)s^{n+d-1}\ee
 and 
 \be\theta(z)t^{n+d-1}\le \H^{n+d-1}(B(z,t)\cap E')\le \frac32\theta(z)t^{n+d-1}.\ee

As a result, we know that 
\be \H^{n+d-1}(E'\cap B(z, s)\bs B(z,t))\le \frac32\theta(z)s^{n+d-1}-\theta(z)t^{n+d-1}.\ee
We apply the coarea formula (\cite{Fe} 3.2.22) to the 1-Lipschitz map $f=d(z,\cdot)$, and get that 
\be \int_{r=t}^s \H^{n+d-2}(E'\cap \pa B(z,r))dr\le \frac32\theta(z)s^{n+d-1}-\theta(z)t^{n+d-1}.\ee
When $s=2t$, we have
\be \int_\frac s2^s \H^{n+d-2}(E'\cap \pa B(z,r))dr\le \frac32\theta(z)s^{n+d-1}-\theta(z)(\frac s2)^{n+d-1}\le \frac32Cs^{n+d-1}.\ee
This means, by Chebyshev, that the set $\{r\in [\frac s2,s]: \H^{n+d-2}(E'\cap \pa B(z,r))\le 8Cs^{n+d-2}\}$ is of positive measure. On the other hand, let $[E'_{w'}]_{C^I_{n+d-1}}$ be the $n+d-1$-integral current $\H^{n+d-1}\lf_{E'}\wg w'_*$, then by slicing for integral current (\cite{Fe} 4.2.1), since $B(z,r)\cap spt(\pa [E'_{w'}]_{C^I_{n+d-1}})=\emptyset$, we know that for a.e. $r\in [\frac s2,s]$, the set $E'\cap \pa B(z,r)$ is $n+d-2$-rectifiable, and $spt(\pa([E'_{w'}]_{C^I_{n+d-1}}\lf_{B(z,r)})\subset E'\cap \pa B(z,r)$. As a result, there exists $r\in [\frac s2,s]$ so that 
\be E'\cap \pa B(z,r)\mbox{ is }n+d-2\mbox{-rectifiable, } \H^{n+d-2}(E'\cap \pa B(z,r))\le C_1 r^{n+d-2},\ee
and \be spt(\pa([E'_{w'}]_{C^I_{n+d-1}}\lf_{B(z,r)})\subset E'\cap \pa B(z,r).\ee
where $C_1= 2^{n+d+1}C$ only depends on $n$ and $d$.

 Now since $\H^{n-2}(S_w)=0$, and $S_w\subset E$, we know that $\H^{n+d-2}(S'_w)=\H^{n+d-2}(S_w\times V)=0$. 
 In particular, $\H^{n+d-2}(S'_w\bs \cT_{2\d})=0$. Moreover we know that $S'_w\bs \cT_{2\d}$ is compact. 
 
As a result, for a fixed $\tau>0$, there exists a finite family of balls $B_i=B(z_i,s_i), 1\le i\le N$ in $U\times V$ with $s_i<\d/2$, so that $S'_w\bs \cT_{2\d}\subset\cup_{i=1}^N B_i$, and $\sum_{i=1}^N s_i^{n+d-2}<2^{2-n-d}\tau$.  We can suppose that $B_i\cap S_w\ne \emptyset$, and hence by replacing $s_i$ by $2s_i$ if necessary, we can suppose that $z_i\in S_w$. Now by the above argument, for each $i$ there exists $r_i\in [s_i,2s_i]$ so that \tb{(4.17) and (4.18)} hold for $z=z_i, r=r_i$:
\be E'\cap \pa B(z_i,r_i)\mbox{ is }n+d-2\mbox{-rectifiable, } \H^{n+d-2}(E'\cap \pa B(z_i,r_i))\le C_1 r_i^{n+d-2},\ee
and \be spt(\pa([E'_{w'}]_{C^I_{n+d-1}}\lf_{B(z_i,r_i)})\subset E'\cap \pa B(z_i,r_i).\ee
Moreover since $r_i\le 2s_i$, we have
\be \sum_{i=1}^N r_i^{n+d-2}<\tau.\ee

Let $B_i'=B(z_i, r_i)$, and let $U_0=\cup_{i=1}^N B'_i$. Then it is an open subset of $U\times V\bs\cT_{2\d}$. 
Note that $E'\cap \pa U_0\subset \cup_{i=1}^N (E'\cap \pa B'_i)$, hence by \tb{(4.21)},
\be\H^{n+d-2}(E'\cap \pa U_0)\le \sum_{i=1}^N \H^{n+d-2}(E'\cap \pa B'_i)\le C_1\sum_{i=1}^N  r_i^{n+d-2}<C_1\tau.\ee

Set $\cG_0=p'^{-1}(E'\bs (U_0\cup \cT_{2\d}))$. It is well defined, because $S'_w\bs \cT_{2\d}\subset U_0$. Let $\cG=\cG_0\cup \cT_{3\d}$. Let $R^2_i=R^1_i\lf_\cG$. Since $spt(R^1_i)\cap \cT_{4\d}\subset E'$, we know that $spt(R^2_i)\subset D'$. Recall that $\pa R_i^1=T_i+(-1)^n[\cU_i]_{C^I_{n+d-2}}$, and $spt(T_i)\subset A$ is $n+d-2$-essentially disjoint from $\cU_i\subset \pa (U\times V)$. Also we have $\cU_i\subset \cT_{2\d}$, hence $\cU_i\cap \pa\cG=\emptyset$. Therefore we have 
\be \cU_i\subset spt(\pa R^2_i)\subset [\cU_i\cup spt(T_i)\cup (\pa\cG\cap spt(R^2_i))]\cap D'.\ee

Let $Q_i=p'_\s(R^2_i)$. Then it is supported in $E'$, and 
\be\begin{split} spt(\pa Q_i)&=spt(\pa p'_\s(R^2_i))=spt(p'_\s\circ\pa R^2_i)\subset p'(spt(\pa R^2_i))\\
&\subset p'([\cU_i\cup spt(T_i)\cup (\pa\cG\cap spt(R^2_i))]\cap D')\\
&=p'(\cU_i)\cup p'(spt(T_i)\cap D')\cup p' (\pa\cG\cap spt(R^2_i))\\
&\subset \cU_i\cup p'(A\cap D')\cup p' (\pa\cG\cap spt(R^2_i)).\end{split}\ee

Note that $\pa \cG=\pa (\cG_0\cup \cT_{3\d})\subset [\pa\cG_0\bs \cT_{3\d}] \cup [\pa \cT_{3\d}\bs \cG_0]$, hence
\be p'[\pa\cG\cap spt(R^2_i)]\subset p'[spt(R^2_i)\cap \pa\cG_0\bs \cT_{3\d}]\cup p'[spt(R^2_i)\cap \pa \cT_{3\d}\bs \cG_0].\ee 

For the first term of \tb{(4.25)}, since $\pa\cG_0\subset U\times V$, we have 
\be\begin{split} p'[spt(R^2_i)\cap \pa\cG_0\bs \cT_{3\d}]&\subset p'[spt(R^2_i)\cap \pa\cG_0\cap (U\times V)]\\
&\subset p'[\pa\cG_0\cap (U\times V)]
\subset (E'\cap \pa U_0);
\end{split}\ee 
For the second term of \tb{(4.25)}, we know that $spt(R^2_i)\cap \cT_{4\d}\subset E'\cap \cG$, and $E'\cap\pa \cT_{3\d}\cap \cG\subset \cG_0$, hence $spt(R^2_i)\cap \pa \cT_{3\d}\bs \cG_0=\emptyset$, and thus
\be p'[spt(R^2_i)\cap \pa \cT_{3\d}\bs \cG_0]=\emptyset.\ee

Combining \tb{(4.24)-(4.27)} we get
\be spt(\pa Q_i)\subset \cU_i\cup p'(A\cap D')\cup (E'\cap \pa U_0).\ee
As a result, by \tb{(4.22)}, 
\be \H^{n+d-2}(spt(\pa Q_i)\bs [p'(A\cap D')\cup \cU_i])\le \H^{n+d-2}(E'\cap \pa U_0)<C_1\tau.\ee
Let $\tau=\e/C_1$, and we get the conclusion.\qed

\begin{cor}Under the hypotheses of Theorem \tb{4.1}. Let $A\in \F_{ihc}^{IR,d}(E\times F, n+d-2,\Gamma,H,U\times V)\}$ so that there exists $\d>0$ such that $A\cap \cT_{4\d}\subset E'$. Then for any $\e>0$, and any $1\le i\le k$, there exists a $n+d-1$ integral current $Q_i$ supported on $E'$, and a $n+d-2$-integral current $K_i$ in $E'$
whose support is $\H^{n+d-2}$ essentially disjoint from $\pa(U\times V)$, so that $\pa Q_i=K_i+[\cU_i]_{C^I_{n+d-2}}$, and $\H^{n+d-2}(p'(A\cap D')\bs K_i)<\e$.
\end{cor}

\nd This is a direct corollary from the above proposition. Just let $K_i=\pa Q_i-[\cU_i]_{C^I_{n+d-2}}$.\qed
 
\noindent\textbf{Proof of Theorem \tb{4.1}.}

After Proposition \tb{1.19} and Proposition \tb{4.2}, it is enough to prove that for each $A\in \F_{ihc}^{IR,n+d-2}(E\times F, n+d-2,\Gamma,H,U\times V)$ so that there exists $\d>0$ such that $A\cap \cT_{4\d}\subset E'$, we have $\H^{n+d-2}(E\times F)\le \H^{n+d-2}(A)$.

So take a such $A$. Let $\e>0$ be given. By Corollary \tb{4.4}, for each $1\le i\le k$, there exists a $n+d-1$ integral current $Q_i$ supported on $E'$, and a $n+d-2$-integral current $K_i$ supported in $\bar U\times \bar V$
whose support is $\H^{n+d-2}$ essentially disjoint from $\pa (U\times V)$, so that $\pa Q_i=K_i+[\cU_i]_{C^I_{n+d-2}}$, and $\H^{n+d-2}(p'(A\cap D')\bs K_i)<\e$. Without loss of generality, we can suppose that $spt(Q_i)$ is connected. Otherwise, we just replace it by its restriction to the connected component that contains $\cU_i$.


Note that the $Q_i, 1\le i\le k$ are of full dimension in $E'$, $E'\cap (U\times V)$ is a manifold except at a set of codimenion at most 8, and the multiplicity of $\pa Q_i$ at $\cU_i$ is 1, so we know that the multiplicity of $Q_i$ near $\cU_i$ is 1. Then by the constancy theorem, since $spt(Q_i)$ is connected, we know that the multiplicity of $Q_i$ is 1 everywhere on $spt(Q_i)$. So let $Q_i=\H^{n+d-1}\lf_{spt(Q_i)}\wg \zeta_i$, where $\zeta_i$ is a measurable orientation on $spt(Q_i)$. Then again since $\pa Q_i=[\cU_i]_{C^I_{n+d-2}}+K_i$, we know that $\zeta_i=w'_*$, $\H^{n+d-1}\lf_{spt(Q_i)}$-a.e.

Set $R_1=Q_1$; and for $2\le i\le k-1$, set $R_i=\H^{n+d-1}\lf_{spt(Q_i)\bs\cup_{j=1}^{i-1}spt(Q_j)}\wg w'_*$; finally set $R_k=\H^{n+d-1}\lf_{(E\times\bar V)\bs\cup_{j=1}^{k-1}spt(Q_j)}\wg w'_*$.

Then the $spt(R_i),1\le i\le k$ are $\H^{n+d-1}$-essentially disjoint, $\cup_{i=1}^k spt(R_i)=E'$, and
\be spt(\pa R_i)\subset \cup_{j=1}^i \spt(\pa R_j)\subset [\cup_{i=1}^k spt(K_i)]\cup [\cup_{i=1}^k \cU_i],\ee
Moreover, since the set $E\times \bar\O_i$ is $\H^{n+d-2}$ essentially disjoint from the $E\times \bar\O_j$ for $j\ne i$, and from the parts $[\cup_{i=1}^k spt(K_i)]$, therefore by definition of $R_i$, we know that last part in the right-hand-side of \tb{(4.30)} is in fact only $spt(\cU_i)$.  That is,
\be\begin{split} spt(\pa R_i)\subset[\cup_{i=1}^k spt(K_i)]\cup \cU_i
\end{split}\ee
where the last union is a $\H^{n+d-2}$-essentially disjoint union. Set $K_i'=\pa R_i-[\cU_i]_{C^I_{n+d-2}}$. Then $spt(K_i')\subset [\cup_{i=1}^k spt(K_i)]$, and $spt(K_i')$ is $\H^{n+d-2}$-essentially disjoint from $\cU_i$.

Since we have the essentially disjoint union $E'=\cup_{i=1}^k spt(R_i)$, and since $E'\cap (U\times V)$ is a manifold except a subset of dimension at most $n-8$, we know that for $\H^{n+d-2}$-a.e. $z\in \cup_{i=1}^n spt(\pa R_i)\cap [E'\cap (U\times V)]=\cup_{i=1}^n spt(K_i')$, there exists exactly 2 $1\le i<j\le k$, so that $z\in spt(\pa R_i)\cap spt(\pa R_j)$.

%

As a result, we have

(1) Let $A'=\cup_{i=1}^k spt(K_i')\subset E'$. Then for $\H^{n+d-2}$-a.e. $z\in A'$, there exists exactly 2 $1\le i<j\le k$, so that $z\in A'_{ij}:=spt(K_i')\cap spt(K_j')$, and $\sigma'_i(z)=-\sigma'_j(z)$, where $\sigma'_i$ is the orientation of $K_i'$. Hence $A'$ is the $\H^{n+d-2}$-essentially disjoint unions $A'=\cup_{1\le i<j\le k}A'_{ij}$, and $spt(K_i')=\cup_{j\ne i}A'_{ij}$.

(2) Moreover, by definition of $K_i'$ and $A'$, we know that $\H^{n+d-2}(A'\bs p'(A\cap D'))<\e$, and hence we know that $\H^{n+d-2}(A'_{ij}\bs p'(A\cap D'))<\e$ for any $1\le i<j\le k$.

Now since the $n+d-2$ forms $w\wg v_i^\perp$ are closed, we have for each $1\le i\le k$,
\be 0=\lg\pa R_i, w\wg v_i^\perp\rg=\lg K_i', w\wg v_i^\perp\rg+\lg [\cU_i]_{C^I_{n+d-2}}, w\wg v_i^\perp\rg.\ee
That is,
\be \lg [\cU_i]_{C^I_{n+d-2}}, w\wg v_i^\perp\rg=-\int_{spt(K_i')}\lg \sigma'_i, w\wg v_i^\perp\rg d\H^{n+d-2}.\ee
By (1), we have
\be \lg [\cU_i]_{C^I_{n+d-2}}, w\wg v_i^\perp\rg=-\sum_{j\ne i}\int_{A'_{ij}}\lg \sigma'_i, w\wg v_i^\perp\rg d\H^{n+d-2}.\ee
We sum over $1\le i\le k$, and since $\sigma'_i=-\sigma'_j$ for $\H^{n+d-2}$-a.e. on $A'_{ij}$, we get
\be \begin{split}&\sum_{1\le i\le k}\lg [\cU_i]_{C^I_{n+d-2}}, w\wg v_i^\perp\rg=-\sum_{1\le i\le k}\sum_{j\ne i}\int_{A'_{ij}}\lg \sigma'_i, w\wg v_i^\perp\rg d\H^{n+d-2}\\
=&\sum_{1\le i<j\le k}\int_{A'_{ij}}\lg \sigma'_i, w\wg v_j^\perp-w\wg v_i^\perp\rg d\H^{n+d-2}.
\end{split}\ee

Now let $A_{ij}$ be the first inverse image of $A'_{ij}\cap p'(A\cap D')$ under $p'$. Then since the $A'_{ij}$ are pairwise disjoint, so are $A_{ij}$. Moreover by (2), we know that
\be\H^{n+d-2}(A'_{ij}\bs p'(A_{ij}\cap D'))<\e, 1\le i<j\le k.\ee

Now for each $i$, let $\sigma_i$ be the orientation on $A_{ij}$, so that $\sigma_i'$ is the orientations on $p(A_{ij})$ associated to $\sigma_i$, as defined in Lemma \tb{2.6}. Then for each $1\le i<j\le k$, by Lemma \tb{2.6 $4^\circ$}, with $\varphi=p'$, $u'=\sigma'_i$ (resp. $\sigma'_j$) on $p'(A_{ij}\cap D')\subset A_{ij}'$, $u=\sigma_i$ (resp. $\sigma_j$), we have

\be \int_{p'(A_{ij}\cap D')}\lg \sigma'_i, w\wg v_i^\perp\rg d\H^{n+d-2}=\int_{A_{ij}}\lg \sigma_i, w\wg v_i^\perp\rg d\H^{n+d-2}\ee
and 
\be \int_{p'(A_{ij}\cap D')}\lg \sigma_j', w\wg v_i^\perp\rg d\H^{n+d-2}=\int_{A_{ij}}\lg \sigma_j, w\wg v_i^\perp\rg d\H^{n+d-2}.\ee
Since $\sigma_i'=-\sigma_j'$ for $\H^{n+d-2}$-a.e. $z\in A'_{ij}$, we have
\be \int_{p'(A_{ij}\cap D')}\lg \sigma_i', w\wg (v_j^\perp-v_i^\perp)\rg d\H^{n+d-2}=\int_{A_{ij}}\lg \sigma_i, w\wg (v_j^\perp-v_i^\perp)\rg d\H^{n+d-2}.\ee

Combine with \tb{(4.35)}, we have
\be \begin{split}&\sum_{1\le i\le k}\lg \cU_i, w\wg v_i^\perp\rg
=\sum_{1\le i<j\le k}\int_{A'_{ij}}\lg \sigma'_i, w\wg v_j^\perp-w\wg v_i^\perp\rg d\H^{n+d-2}\\
=&\sum_{1\le i<j\le k}[\int_{p'(A_{ij}\cap D')}\lg \sigma'_i, w\wg (v_j^\perp-v_i^\perp)\rg d\H^{n+d-2}+\int_{A'_{ij}\bs p'(A_{ij}\cap D')}\lg \sigma'_i, w\wg (v_j^\perp-v_i^\perp)\rg d\H^{n+d-2}]\\
=&\sum_{1\le i<j\le k}\int_{A_{ij}}\lg \sigma_i, w\wg (v_j^\perp-v_i^\perp)\rg d\H^{n+d-2}+\sum_{1\le i<j\le k}\int_{A'_{ij}\bs p'(A_{ij}\cap D')}\lg \sigma'_i, w\wg (v_j^\perp-v_i^\perp)\rg d\H^{n+d-2}\\
\le &\sum_{1\le i<j\le k}\H^{n+d-2}(A_{ij})||w\wg (v_j^\perp-v_i^\perp)||_\infty d\H^{n+d-2}+\sum_{1\le i<j\le k} \H^{n+d-2}(A'_{ij}\bs p'(A_{ij}\cap D'))||w\wg (v_j^\perp-v_i^\perp)||_\infty\\
\le &\sum_{1\le i<j\le k}M(\nu)\H^{n+d-2}(A_{ij})+\tbinom{k}{2}M(\nu)\e\le M(\nu)\H^{n+d-2}(A)+\tbinom{k}{2}M(\nu)\e
\end{split}\ee
by \tb{(4.36)}, and since the $A_{ij}, 1\le i<j\le k$ are disjoint subsets of $A$.

On the other hand, if we take $A=E\times F$, then it is easy to see that $A_{ij}=A'_{ij}=E\times F_{ij}$, and hence
\be \sum_{1\le i\le k}\lg [\cU_i]_{C^I_{n+d-2}}, w\wg v_i^\perp\rg=M(\nu)\H^{n+d-2}(E\times F)\ee
since $F\in PC(\mo, \nu)$. Hence we have, for any $A\in \F_{ihc}^{IR,n+d-2}(E\times F, n+d-2,\Gamma,H,U\times V)$ so that there exists $\d>0$ such that $A\cap \cT_{4\d}\subset E'$, and every $\e>0$, we have
\be M(\nu)\H^{n+d-2}(E\times F)\le M(\nu)\H^{n+d-2}(A)+\tbinom{k}{2}M(\nu)\e.\ee
Thus 
\be \H^{n+d-2}(E\times F)\le \H^{n+d-2}(A).\ee
This completes our proof of Theorem \tb{4.1}.\qed

\begin{thm} Let $n\ge 2,d\ge 2$ be integers. 

Let $n\ge 2,d\ge 2$ be integers. 

Let $U$ be an open bounded set in $\R^n$. Let $\b$ be an $n-2$-integral current without boundary, such that $B=spt(\b)\subset \pa U$.  Let $E\in C_1(\b,w,U)$, where $w$ is a coflat calibration in $U$, such that its singular set $S_w\subset E$, and $\H^{n-2}(S_w)=0$. 

Let $V\subset \R^d$ be a bounded strongly Lipschitz domain, so that $(\R^d,\pa V)$ admits a $C^1$ triangulation. Let $F\in PC(\mo,\nu)$, where $\mo=(\Omega_1,\cdots, \Omega_k)\in \cO_P(V,k)$ for some $k\in \N$, and $\nu\in \cV_c^k(\bar V)$ is such that $v_i-v_j$ does not vanish on $V\bs S_\nu$ for any $i<j$.

Suppose that $U\times V$ is a SR-domain. 

Then 

$1^\circ$ The product $E\times F$ is a $n+d-2$-dimensional Almgren minimal set in $U\times V$.

$2^\circ$ Suppose in addition that the $(n+d-2)$-th integral rectifiable homology group $H^I_{n+d-2}(\bar U\times\bar V)$ is trivial.
Let $H=\{\pa[[E_w]_{C^I_{n-1}}\times [F_i]_{C^I_{d-1}}], 1\le i\le k\}$, where $F_i$ and $[F_i]_{C^I_{d-1}}$ are as defined at the beginning of the section.
Then  \be \H^{n+d-2}(E\times F)=\inf\{\H^{n+d-2}(A):A\in \F_{isc}(n+d-2,(B\times F)\cup (E\times B(\mo)),H,U\times V)\}.\ee
\end{thm}

\nd $1^\circ$ is a direct corollary of Theorem \tb{4.1} and Proposition \tb{1.16}; $2^\circ$ follows directly from Theorem \tb{4.1}.\qed

\section{New types of singularities for Plateau's problem}

As stated in the introduction, we can find new types of Almgren minimal cones with the help of Theorem \tb{4.5}. We will give examples below. Note that for any cone, its minimality is equivalent to its minimality in any open convex domain that contains the origin. And the regularity of any convex domain is quite enough for us to apply our theorems in the previous sections.

Let us first introduce some notations:

-A set $C\subset \R^n$ is called a cone, if it contains $0$, and for any $\lambda>0$, we have $\lambda C=C$. In other words, it is a set composed of rays issued from the origin.

-Given a cone $C\subset \R^n$, it is completely determined by its intersection with the unit sphere $X=C\cap \mS^{n-1}$. 

-Conversely, given any subset $X\subset  \R^n$, we let $C(X)=\{tx: 0\le t<\infty, x\in X\}$ denote the cone over $X$. In particular, if $X\subset \mS^{n-1}$, then $C(X)\cap  \mS^{n-1}=X$.

-If $X\subset \mS^{n-1}$ is an oriented $n-2$-submanifold, then $C(X)$ is called a regular cone. It is a manifold except at the origin.

\subsection{Examples of coflat calibrated minimal cones of multiplicity 1}

\begin{pro}[Multiplicity for regular calibrated area minimizing hypercones]Let $M\subset \mS^{n-1}$ be an oriented $n-2$-submanifold. Let $C(M)$ be the cone over $M$. If $w$ is a coflat calibration in the unit ball $B^n$, so that $C(M)$ is calibrated by $w$ with respect to $[M]_{C^I_{n-2}}$, then $C(M)$ is of multiplicity 1. In particular, $C(M)\in C_1([M]_{C^I_{n-2}}, w,B)$.\end{pro}

\nd By Proposition \tb{2.4}, we know that $C(M)$ is of multiplicity no less than 1. Hence it is enough to prove that the multiplicity is no more than one. 

Since $M$ is a closed surface of codimension 1, we know that $C(M)$ separates the unit ball $B^n$ into two connected components $C_1$ and $C_2$. Let $u$ be the orientation of $C(M)$ induced by the orientation of $M$. Without loss of generality, we can suppose that the vector field $u_\perp$ points ourward from $C_1$ into $C_2$. Again by Proposition \tb{2.4}, we know that $\lg u_\perp,w_\perp\rg=1$ $\H^{n-1}$-a.e. on $C(M)$. Now for any integral curve $\gamma$ of $w_\perp$ in $B\bs S_w$, if it touches $C(M)$ at a point $p$, then since $\lg u_\perp,w_\perp\rg=1$, we know that it goes from $C_1$ to $C_2$ via $p$. As a result, it cannot touch $C(M)$ at more than one point because if so, then it goes from $C_1$ to $C_2$ more than once without going from $C_2$ to $C_1$, which is not possible.\qed

As a corollary of the above proposition, we have the following coflat calibrated minimal cone of multiplicity 1 that satisfy the condition of Theorems \tb{4.1 and 4.5}:

\begin{exam}[Regular area-minimizing hypercones]
$1^\circ$ Every area-minimizing cone $C_{r-1, s-1}\subset \R^{r+s}$ (in particular, the Simons cone) (cf. \cite{ZYS16} Theorem 1.8);

$2^\circ$ More generally homogeneous area-minimizing hypercones $C$ (cf. \cite{ZYS16} Theorem 1.9).
\end{exam}

\subsection{Examples of paired calibrated minimal cones}

In this subsection we simply give examples of coflat paired calibrated cones.

\begin{exam}[Regular calibrated hypercones of multiplicity 1]As in the proof of Proposition \tb{5.1}, we know that every regular hypercone $C(M)$ separates the unit sphere $B$ into two connected components $\O_1,\O_2$, where $M$ is an oriented $n-2$ submanifold of $\mS^{n-1}$. Then $C(M)$ separates $B$ into two connected components $C_i=C(\O_1)\bs \{0\}, i=1,2$. Now suppose $w$ is a coflat calibration such that $C(M)\in C([M]_{C^I_{n-2}}, w,B)$. Modulo changing the index we can suppose that on almost every point on $M$, $w_\perp$ points ourward from $C_1$ into $C_2$. By the definition of a coflat calibration and Proposition \tb{2.4}, we know that $||w||_\infty=1$.

Now $\mo=(\O_1,\O_2)\in \cO_P(B, 2)$ is a Lipschitz partition of $\pa B$, $\nu=(w_\perp, -w_\perp)\in \cV_c^2(\bar B)$, and $M(\nu)=||w_\perp-(-w_\perp)||_\infty=||2w_\perp||_\infty=2$. Then $C(M)\in PC(\mo,\nu)$. 
\end{exam}

\begin{exam}[Cones over simplices in $\R^n$ with $d\ge 3$, \cite{LM94}] Let $d\ge 3$. Let $\sigma$ be a unit regular $d$-simplex centered at the origin of $\R^d$. Let $C_\sigma$ be the cone over the $d-2$-skeleton of $\sigma$. Let $\mo$ be the Lipschitz partition of $\pa \sigma$ so that $B(\mo)=C_\sigma\cap \pa \sigma$. Then there exists a family of coflat paired calibrations $\nu$ so that $C_\sigma\in PC(\mo,\nu)$.
\end{exam}

\begin{exam}[Cones over squares in $\R^d$ with $d\ge 4$, \cite{Br91}] Let $d\ge 4$. Let $Q$ be the closed unit cube centered at the origin of $\R^d$. Let $C_Q$ be the cone over the $d-2$-skeleton of $Q$. Let $\mo$ be the Lipschitz partition of $\pa Q$ so that $B(\mo)=C_Q\cap \pa Q$. Then there exists a family of coflat paired calibrations $\nu$ so that $C_Q\in PC(\mo,\nu)$.
\end{exam}

\subsection{Examples of new families of singularities}

Take the examples in the previous 2 subsections and apply Theorem \tb{4.5}, we have:

\begin{thm}

$1^\circ$ The product of any two homogeneous area-minimizing hypercones (Example \tb{5.2}) (same or different) is an Almgren minimal cone. 

$2^\circ$ The product of any homogeneous area-minimizing hypercone (Example \tb{5.2, 5.3}) and any cone over a simplex in $\R^d$ for $d\ge 3$ (Example \tb{5.4}) is an Almgren minimal cone.

$3^\circ$ The product of any homogeneous area-minimizing hypercone (See Example \tb{5.2, 5.3}) and any cone over a cube in $\R^d$ for $d\ge 4$ (Example \tb{5.5}) is an Almgren minimal cone.
\end{thm}

\renewcommand\refname{References}
\bibliographystyle{plain}
\bibliography{reference}

\end{document}